\documentclass[12pt]{article}

\usepackage{amsmath}
\usepackage{amssymb}
\usepackage{amsfonts}
\usepackage{amsthm}
\usepackage{enumerate}

\newtheorem{theorem}{Theorem}[section]
\newtheorem{proposition}[theorem]{Proposition}
\newtheorem{lemma}[theorem]{Lemma}
\newtheorem{definition}[theorem]{Definition}
\newtheorem{corollary}[theorem]{Corollary}

\newcommand{\bref}[1]{(\ref{#1})}

\renewcommand{\l}{\lambda}
\newcommand{\e}{\epsilon}

\newcommand{\M}{{\cal M}}
\renewcommand{\P}{{\cal P}}

\renewcommand{\d}{\delta}

\newcommand{\R}{{\mathbb R}}
\newcommand{\rr}{{\mathbb R}}
\newcommand{\Z}{{\mathbb Z}}

\newcommand{\T}{{\cal T}}

\newcommand{\calr}{{\cal R}}

\newcommand{\lo}[1]{{(\log\frac{1}{\d})^{#1}}}

\newcommand{\supp}{\hbox{supp\,}}
\newcommand{\sdel}{S_\delta}

\newcommand{\ssd}{\Sigma^\sigma_\d}
\linethickness{.4mm}

\begin{document}

\title{Wolff's inequality for hypersurfaces}

\author{Izabella {\L}aba and Malabika Pramanik}

\maketitle 

\begin{abstract}
We extend Wolff's ``local smoothing" inequality \cite{Wls} to
a wider class of not necessarily conical hypersurfaces of 
codimension 1.  This class includes surfaces with nonvanishing
curvature, as well as certain surfaces with more than one
flat direction.  An immediate consequence is the $L^p$-boundedness
of the corresponding Fourier multiplier operators.

\end{abstract}

The purpose of this article is to extend the ``local smoothing"
inequality, proved in \cite{Wls}, \cite{LW} for circular cones in $\R^d$,
$d\geq 2$, 
and in \cite{PS} for more general conical surfaces in $\R^3$, 
to a wider class of bounded surfaces of codimension $1$ in 
$\R^{d+1}$, $d\geq 3$.  

Recall that Wolff's inequality \cite{Wls} states that if $f$ is a function 
with Fourier transform supported in a $\d$-neighbourhood of the
segment of the circular cone with $1\leq |x|\leq 2$, then
\begin{equation}\label{cone-est}
\forall\e\,\exists C_{\e}:\;\|f\|_p\leq C_{\e}\d^{-\frac{d-1}{2}+\frac{d}{p}-\e}
\|f\|_{p, \d},
\end{equation}
with $d=2$ and $p>74$; this was then extended in \cite{LW} to
$d\geq 3$, $p>2+\frac{32}{3d-7}$, and
$d\geq 4$, $p>2+\frac{8}{d-3}$.
Here the norm on the right side is defined by
\[
\|f\|_{p, \d}=\left(\sum_{a}\|\Xi_a\ast f\|_p^p\right)^{1/p},
\]
where $\widehat{\Xi_a}$ are, roughly, cutoffs in the Fourier space corresponding
to the natural covering of the $\delta$-neighbourhood of the cone by 
rectangular ``plates" of thickness $\d$.  We propose to extend this
inequality in two directions.
First, we will consider cones generated by more general surfaces 
of codimension $2$.  Second, we will also allow surfaces with more than
one flat direction, satisfying certain geometrical conditions
stated below; the simplest nontrivial examples are the ``$k$-cones",
defined below and described in more detail in Section \ref{k-cones}.
An immediate consequence is the $L^p$ boundedness of Bochner-Riesz
multipliers for the same classes of surfaces and for appropriate
ranges of exponents.  

We now define specific classes of surfaces of interest to us.  
We will always assume that $S$ is a bounded surface of codimension
$1$ in $\rr^{d+1}$, smooth everywhere except for
the possible boundary, with all curvatures bounded from
above uniformly by a constant.

A {\em nondegenerate surface} in $\R^{l+1}$ will be a surface 
$S_0$ defined by an equation of the type $x_{l+1}=F(x_1,\dots,x_l)$, 
where $F$ is smooth with all derivatives bounded uniformly in $x$,
and all principal curvatures are bounded away from zero:
$$|\langle u, (D^2_x F)u\rangle| \geq c_0|u|^2,\ u\in \R^{l+1},$$
with the lower bound $c_0$ uniform in $u$.  
Elliptic surfaces, such as e.g. a sphere or a paraboloid, are clearly 
nondegenerate.  However, we do not require the Hessian to be positive 
definite, hence we also allow surfaces with both positive and
negative curvatures (e.g. hyperboloids).

A {\em conical surface} in $\R^{d+1}$ is defined as 
\begin{equation}\label{def-cone}
S=\{tx: \ x\in S_0,\ t\in [C_1,C_2]\},\ 0<C_1<C_2<\infty,
\end{equation}
where $S_0$ is a surface of dimension $l=d-1$ 
contained in an affine subspace $X\subset \rr^{d+1}$ of dimension
$d$ which does not pass through the origin, such that 
$S_0$ viewed as a subset of $X$ is a nondegenerate surface
if $X$ is identified with $\rr^d$ in the obvious way.
This class includes circular cones as well as more general 
homogeneous quadrics of the form 
$S=\{x:\ 1\leq |x|\leq 2,\ \langle Ax,x\rangle =0\}$,
where $A$ is a symmetric $(d+1)\times (d+1)$ matrix of full rank;
note that we {\em do not} assume anything about the signature of $A$.

A {\em $k$-cone} in $\R^{d+1}$, where $1\leq k\leq d-1$,
is constructed as follows.   Let $L_0$
be a $(d-k+1)$-dimensional linear subspace of $\rr^{d+1}$, and let 
$L_i=L_0+v_i$ for $i=1,\dots,k$, where $v_1,\dots,v_{k}$ are
linearly independent vectors such that $L_0,v_1,\dots,v_k$ span
$\rr^{d+1}$.  In each of the subspaces $L_i$ we fix a bounded
strictly convex solid $F_i$ such that $E_i=\partial F_i$ is smooth and
has nonvanishing Gaussian curvature.  Thus $E_i$ is a surface
of dimension $d-k$ (i.e. of codimension $1$ in $L_i$).  
We say that a $(k+1)$-tuple of points $(x_0,x_1,\dots,x_k)$
is {\em good} if $x_i\in E_i$, $i=0,\dots,k$, and if the outward
unit normal vectors to $E_i$ at $x_i$ are the same.  We 
then let
$$
S=\bigcup_{(x_0,\dots,x_k)\ good}\eta(x_0,\dots,x_k),
$$
where $\eta(x_0,\dots,x_k)$ denotes the convex hull of $x_0,\dots,
x_k$ in $\rr^{d+1}$.  In Section \ref{k-cones} we prove that this indeed defines
a surface of codimension 1 and that each point $a\in S$ 
belongs to $\eta(x_0,\dots,x_k)$ for exactly one $(k+1)$-tuple
$(x_0,\dots,x_k)$; we will then call $\eta(x_0,\dots,x_k)$
the {\em $k$-plane } at $a$, and denote it by $\eta(a)$.

For illustration purposes, consider the simple case when
$k=1$ and $L_0$ and $L_1$ are two
parallel hyperplanes.  If $E_0$, $E_1$ are spheres of 
different radii, then $S$ is a segment of a right circular cone
or a slanted circular cone, depending on the relative location of $E_0,E_1$.  
Similarly, if $E_0,E_1$ are spheres of equal radii, $S$ is a segment of a 
(right or slanted) circular cylinder.  However, if $E_0,E_1$ are
randomly chosen ellipsoids, then $S$ will usually not be a cone,
a cylinder, or an affine image thereof.  For $k\geq 2$, these surfaces
may be more difficult to visualize; see Section \ref{k-cones} for
more details.  It is likely that similar surfaces
may be constructed if $E_i$ are allowed to be more general surfaces of
codimension $1$ in $L_i$, but since it seems difficult to find
the precise conditions on such more general $E_i$ under which the
construction works, we choose not to do this here.  

What we will actually use in the statement and proof of our main
inequality is that these surfaces admit a uniform
``plate covering", analogous to that of \cite{Wls}, \cite{LW}.
Let $S$ be a conical surface or a $k$-cone, as defined above.
For $a\in S$, we use $n(a)$ to denote the unit normal to $S$ at $a$.
We may assume that the map $x\to n(x)$ is continuous; if $S$ is not
orientable, we restrict our attention to an orientable subset of $S$.
We always use $k$ to denote the number of ``flat" directions of $S$.
For $\delta>0$, let $\sdel$ denote the $\delta$-neighbourhood of $S$.
We write $A\lesssim B$ if $A\leq c B$ for some constant $c$ 
independent of $\delta$, and $A\approx B$ if $A\lesssim B$ and
$B\lesssim A$.
Then $S$ satisfies the following conditions.

\medskip

\noindent{\bf Assumption (A)}.
{\em For each $\delta>0$,  $\sdel$ admits a covering 
$\sdel \subset \bigcup_{a\in\M_{\delta}} \Pi_{a,\delta},$
where $\M_{\delta}\subset S$, and $\Pi_{a,\delta}$ are finitely 
overlapping rectangular boxes centered at $a$ with the
following properties:

\begin{itemize}

\item $c\Pi_{a,\delta}\subset S_{\delta}\cap
\{x\in \rr^{d+1}:\ (x-a)\cdot n(a)\leq \delta\}$, for some small $c$.
\item All $\Pi_{a,\delta}$ have dimensions
$C\d\times C\d^{1/2}\times\dots\times C\d^{1/2}\times C\times\dots
\times C$, where the short direction is normal to $S$ at $a$, the 
long directions are parallel to the $k$-plane $\eta(a)$ at $a$, 
and the mid-length directions are tangent to $S$
at $a$ but perpendicular to $\eta(a)$. 
\item (angular separation) For any $a\in\M_\d$, there are at most $O(1)$
distinct points $b\in\M_{\delta}$ such that
$|n(a)-n(b)|\geq c\d^{1/2}$.
\item (consistency) If $0<\delta\leq\sigma$ and if $b\in \M_{\delta}\cap
\Pi_{a,\sigma}$ for some $a\in \M_{\sigma}$, then $\Pi_{b,\delta}\subset
C''\Pi_{a,\sigma}$.

\end{itemize}

}

Here and in the sequel, ``finitely overlapping" means the following:
if a family of sets ${\cal S}_\delta$ is given for each $\delta$, 
any $x\in\R^{d+1}$ belongs to at most $K$ sets in ${\cal S}_\delta$.
$K,C,C',\dots$ denote 
constants independent of $x,a,b,\delta,\sigma$, and the choice of 
$\M_{\delta}$ and $\M_\sigma$.

We call $\Pi_{a,\delta}$ {\em $\d$-sectors} (note that our terminology is slightly 
different from that of \cite{Wls}, \cite{LW}).  
If $S$ is a nondegenerate hypersurface or a conical surface generated by
$S_0$,  we may take $\M_{\delta}$ to be a maximal 
$\d^{1/2}$-separated subset of $S$ or $S_0$, respectively;
the case of $k$-cones is discussed in Section \ref{k-cones}.
Whenever the choice of the small parameter $\delta$
is clear from the context, we will write $\Pi_{a}$ 
instead of $\Pi_{a,\delta}$.
Note that 
\begin{equation}\label{eq-Pi}
|\Pi_{a,\d}|\approx \d^{\frac{d-k}{2}+1}.
\end{equation}

Let $\Xi_a$ be smooth functions such that
$\|\Xi_a\|_1\approx 1$
and $\{\widehat{\Xi_a}\}_{a\in \M_a}$ is a smooth partition
of unity on $\sdel$ with $\supp \widehat{\Xi_a}\subset \Pi_a$.
Note that the latter condition implies that $\widehat{\Xi_a}\approx 1$
on a box $c\Pi_a$ of size about $\d^{\frac{d-k}{2}+1}$, hence
$\|\Xi_a\|_2^2=\|\widehat{\Xi_a}\|_2^2\approx\d^{\frac{d-k}{2}+1}$.
We may thus choose $\Xi_a$ to have size approximately $\d^{\frac{d-k}{2}+1}$ 
on a box dual to $\Pi_a$ of volume
about $\d^{-\frac{d-k}{2}-1}$, so that the $L^1$ estimate is satisfied.

If $\supp\hat{f}\subset\sdel$, we define
\[
\|f\|_{p, \d}=\left(\sum_{a\in\M_a}\|\Xi_a\ast f\|_p^p\right)^{1/p}
\]
for $2\leq p<\infty$, and
\[
\|f\|_{\infty, \d}=\sup_{a\in\M_a}\|\Xi_a\ast f\|_{\infty}.
\]

Our main result is the following: 

\begin{theorem}\label{thm1}
Assume that (A) holds. 
Then for all functions $f$ with $\supp\hat{f}\subset\sdel$
we have the estimate
\begin{equation}
\forall\e\,\exists C_{\e}:\;\|f\|_p\leq C_{\e}\d^{-\frac{d-k}{2}+\frac{d-k+1}{p}-\e}
\|f\|_{p, \d},
\label{qa3}\end{equation}
with $C_\e$ depending only on $\e$ and on the implicit constants in (A), 
provided that $d$, $k$, $p$ satisfy at least one of the following:

\smallskip

(i) $k<d/3$, $p>p_1(d,k):=2+\frac{8}{d-3k}$,

\smallskip
(ii) $k<\frac{3d-3}{4}$, $p>p_2(d,k):=2+\frac{32}{3d-4k-3}$.

\end{theorem}

For the special case of the spherical cone in $\R^{d+1}$, this is the result of \cite{LW}.
For nondegenerate conical surfaces with $d=2$ and $k=1$ (i.e.,
cones in $\mathbb R^3$ generated by plane curves of nonvanishing
curvature) the inequality \bref{qa3} with $p>74$ (the same exponent as in
\cite{Wls}) has already been obtained by Pramanik and
Seeger in \cite{PS}. The technique in their paper was to approximate the cone in
question by circular cones to which a variant of the result of \cite{Wls} could
be applied. Therefore, the main new cases of interest in the present paper
are $k$-cones with $k \geq 2$ and conical surfaces with both positive and
negative curvatures; such surfaces exist in $\R^{d+1}$ only if $d
\geq 3$. 

Observe that (ii) gives a better range of $k$ for all $d\geq 2$.  The range of
$p$ given in (i) is better than (ii) if and only if $d>8k-3$.
The best possible range of $p$ for which \bref{qa3} could be 
expected to hold is
\begin{equation}\label{best-p}
p\geq 2+\frac{4}{d-k}.
\end{equation}
This may be seen by considering the same example as in \cite{Wls},
page 1238: construct a function $f=\sum_{a\in\M_{\d}}f_a$, where
$\widehat{f_a}$ is supported in a small cube of sidelength $\d$
contained in $\Pi_a$, $|f_a|\leq 1$ on a cube $Q$ of sidelength
$\d^{-1}$, $f_a$ is bounded from below on a smaller proportional
cube $cQ$, and decays rapidly outside of $Q$.  Then
\[
\|f\|_2^2\approx\sum\|f_a\|_2^2\approx M_\d|Q|=\d^{-\frac{d-k}{2}}|Q|,
\]
\[
\|f\|_p^p\gtrsim\|f\|_{L^2(Q)}^p\gtrsim\d^{-\frac{(d-k)p}{4}}|Q|.
\]
Now plug this into the inequality \bref{mnmnmn} with $\sigma
\approx 1$, which will be shown to follow from \bref{qa3}, and
take $\alpha\to 0$.  Using also \bref{eq-M} and comparing the
exponents of $\d$ on both sides, we see that \bref{best-p}
must hold.

If $S$ is a conical surface in $\R^{d+1}$ generated by a nondegenerate
surface $S_0$ in $\R^d$, the exponent of $\d$ in \bref{qa3} is 
$-\frac{d-1}{2} +\frac{d}{p}-\e$ both for $S$ and for $S_0$.  This confirms
the observation of \cite{Wls} that one can deduce \bref{qa3} for the
nondegenerate case from the conical case, for those exponents $p$
for which \bref{qa3} for the conical case is available.  (The idea is
to extend a function supported in a $d$-dimensional neighbourhood
of $S_0$ to a homogeneous degree 0 function defined near $S$, and
then apply \bref{qa3} to $S$.)  However, going through the entire
proof with $k=0$ allowed, rather than using the shortcut just described, yields
a slightly better range of $p$, namely $p>2+\frac{8}{d-1}$ as opposed to
$2+\min(\frac{8}{d-3},\frac{32}{3d-7})$.

An immediate corollary of Theorem \ref{thm1} is the following result
concerning the boundedness of Fourier multipliers associated with
$S$, defined in the usual manner: 
\[
\widehat{T_\alpha f}=m_\alpha \widehat{f},\ 
m_\alpha(\xi)=|\hbox{dist}\,(\xi,S)|^\alpha \phi(\xi),
\]
where $\alpha>0$ and $\phi\in C_0^\infty$ is a suitable smooth cut-off function 
supported in a neighbourhood of $S$.

\begin{corollary}\label{multipliers}
Let $S$ be as in Theorem \ref{thm1}.  Then
$T_\alpha$, defined as above, are bounded on $L^p(\rr^{d+1})$ if
\begin{equation}\label{e-br}
\alpha>(d-k+1)\Big|\frac{1}{2}-\frac{1}{p}\Big|-\frac{1}{2},
\end{equation}
and if one of the following holds:

\smallskip

(i) $k<d/3$ and either $p>p_1(d,k)$ or $1\leq p<\frac{p_1(d,k)}{p_1(d,k)-1}$,

\smallskip
(ii) $k<\frac{3d-3}{4}$ and either $p>p_2(d,k)$ or
$1\leq p<\frac{p_2(d,k)}{p_2(d,k)-1}$,

\smallskip\noindent
where $p_i(d,k)$ are as in Theorem \ref{thm1}.
\end{corollary}

The proof for $p>p_i(d,k)$ is identical to that of Corollary 2(ii) in
\cite{Wls}, therefore we do not reproduce it here.  For 
$p<p_i(d,k)/(p_i(d,k)-1)$, the result follows by duality.

The range of $\alpha$ in \bref{e-br} is sharp for a fixed $p$, see
e.g. \cite{Stein}, pp. 389-390, or \cite{Wls} for a discussion of the
cases $k=0,1$.  On the other hand, the range of $p$ is not sharp,
and in particular it is possible to have $L^p$-boundedness of $T_\alpha$ 
for exponents $p>2$ for which \bref{qa3} fails.  
Indeed, the best possible range of $p$ for
\bref{qa3} with $k=0$ is $p\geq 2+\frac{4}{d}$ (see \bref{best-p}), 
but on the other hand the spherical Bochner-Riesz multipliers with
$\alpha$ as in \bref{e-br} are
known to be bounded for $p>2+\frac{4}{d+1}$ and $p<2-\frac{4}{d+5}$ 
if $d\geq 2$ \cite{Lee}, and the conjectured range is all $1\leq p\leq\infty$.
Moreover, in dimension 2 (i.e. $d=1$, $k=0$) $T_\alpha$ are known to be
bounded for the optimal range of exponents $4/3\leq p\leq 4$, $\alpha>0$,
and $\frac{4}{3+2\alpha}<p<\frac{4}{1-2\alpha}$, $0<\alpha\leq 1/2$
(\cite{CS}, \cite{F}, \cite{Hor2}).

For the case $k=1$, we
recover the result of \cite{LW} for circular cones and extend it
to more general surfaces.  We do not know of any earlier results of 
this type for nontrivial (i.e. not cylindrical) surfaces with $k\geq 2$.  
For further discussion of the existing literature on cone
multipliers, we refer the reader to
e.g. \cite{B}, \cite{MSS}, \cite{Stein}, \cite{TV}, \cite{bil}, \cite{Wls}.
Wolff's inequality \bref{qa3} for cones has also been used to deduce
a variety of other results, including an optimal $L^p$ local smoothing 
result for solutions of the wave equation \cite{Wls} and $L^p$ boundedness
of maximal operators associated with curves in $\rr^3$ \cite{PS}.
We plan to explore the applications of \bref{qa3}
for more general surfaces in a future paper.

The proof of Theorem \ref{thm1} essentially follows the ``induction on scales"
arguments in \cite{Wls}, \cite{LW}, with modifications which we now
describe.  The inductive argument of \cite{Wls}, \cite{LW} (rescaled to
our setting) involves an application of \bref{qa3} on scale $\sqrt{\d}$,
then rescaling each $\sqrt{\d}$-sector to a neighbourhood of the entire cone
via a Lorentz transformation, followed by a second application of
\bref{qa3} on scale $\sqrt{\d}$.  In our more general setting, 
a non-homogeneous scaling cannot be expected to map $S$ to itself,
and in particular Lorentz transformations are usually not available.
Instead, we work directly with the two nested decompositions of $S$.
This will require us to prove \bref{qa3} in somewhat greater generality,
allowing for functions with Fourier supports in sectors $S^\sigma_\delta$,
$\d\leq \sigma\ll 1$.  We do this by first dividing the $\sigma$-sector
in question into $\rho$-sectors, 
where $\rho=\sqrt{\sigma\d}$ is our intermediate scale,
and then subdividing the latter into $\d$-sectors.  Thus \bref{qa3}
for $S^\sigma_\d$ is obtained by combining \bref{qa3} for $S^\sigma_\rho$ 
and $S^\rho_\delta$.  

This of course sounds too good to be true, and
it indeed is: if implemented exactly as described, the above argument
would fail due to the accumulation of the $\d^{-C\e}$-errors arising
at each step of the induction.  This problem is resolved as in \cite{Wls},
\cite{LW}.  Namely, we observe that if the $\rho$ in $S^\sigma_\rho$ 
can be replaced by a slightly bigger scale $\rho^{1-\e_0}$ for some 
fixed $\e_0$, we gain additional factors of $\d^{C\e_0}$ which absorb
the troublesome errors.  We then have to find conditions
under which it is possible to do so.  To this end, we decompose $f$
into standardized ``wave packets", Fourier-localized in $\rho$-sectors
and almost localized spatially in the dual plates.  The inequality
\bref{qa3}, for large $p$, is a statement about the size of set of large values
of $f$.  Fix a tiling of $\R^{d+1}$ by $\rho^{-1+\e_0}$-cubes;
then a combinatorial argument, similar in spirit to Bourgain's 
``bush" argument, shows that if $\lambda$ is sufficiently large relative
to the total number of wave packets, then the sets of wave packets contributing
to the parts of $\{|f|\geq\l\}$ localized in different cubes are
essentially disjoint.  We can now 
adjust the scale by discarding the non-contributing part of each packet
(cf. the ``two ends" reduction of \cite{Wolff-K}).
The rest of the proof is arranged so as to make this step possible,
and in particular this is what determines our range of $p$.

\bigskip

{\bf Acknowledgement:}
This work was partially supported by NSERC grant 22R80520 and by 
NSF grant DMS-0245408.


\section{Notation and preliminaries}\label{prel}


\noindent{\bf General notation:}
\medskip

$S$: a $d$-dimensional bounded connected surface in $\rr^{d+1}$, $C^2$ everywhere
except boundary, with all curvatures bounded, satisfying (A).

$p$: an exponent as in Theorem \ref{thm1} which will remain fixed 
throughout the paper.

$r=1-\frac{2}{p}-\frac{2}{p(d-k)}.$

$\delta,\rho,\sigma$: small dyadic parameters, always satisfying $0<\delta\leq
\rho\leq\sigma\lesssim  1$.

$\e_0, \e_1,\e_2,\e_3$ : small positive numbers with $\e_{j+1}$ much smaller
than $\e_{j}$, depending only on $d,k,p$, to be fixed later.
They will remain fixed through Sections 3--6.

$t$: a dyadic number such that $t\approx (\sigma/\d)^{\e_0}$.  

$C,C_i,C'$, etc.:  constants which 
may depend on the choice of $S$ and $p$, and in particular
on the implicit constants in (A), but will always be 
independent of $\delta$ and of the choice of 
sector decomposition in (A).  They may change from line to line and may be 
adjusted as needed, in particular after each application of Proposition \ref{scales}.

$A\lesssim B$: $A\leq C B$ for some constant $C$.

$A\approx B$: $A\lesssim B$ and $B\lesssim A$.  

$A\lessapprox B$:  $A\lesssim \lo{C}B$ for some constant $C$.

$\chi_E$: the indicator function of the set $E$.

$|E|$: the Lebesgue measure or cardinality of $E$, depending on the context.

A {\em logarithmic fraction of $E$}: a subset of $E$ with measure $\gtrapprox|E|$. 

\medskip

An {\em $l$-cube} is a cube of side length $l$ belonging to a suitable fixed 
$l$-grid on $\R^{d+1}$; thus any two $l$-cubes are either identical or have disjoint
interiors.  If $l$ is fixed, for any $x\in\R^{d+1}$ we let $Q(x)$ be an
$l$-cube such that $x\in Q(x)$.

If $R$ is a rectangular box (e.g. a tube or a plate), we will denote by
$cR$ the box obtained from $R$ by dilating it by a factor of $c$ 
about its center.

If $R_0$ is a rectangular box centered at the origin, the {\em dual box}
to $R_0$ is the rectangular box
\[
R_0^*=\{x\in\rr^{d+1}:\ |x\cdot y|\leq 1\hbox{ for all }y\in R_0\},
\]
where $\cdot$ denotes the usual scalar product in $\rr^{d+1}$.  
We will sometimes say that two boxes $R,R^*$, not centered at the origin,
are dual to each other if and only if their translates $R_0,R_0^*$ 
centered at the origin are dual to each other.

We let $\phi(x)=(1+|x|^2)^{-\frac{K}{2}}$ with $K$ large enough,
and $\phi_R=\phi\circ u_R^{-1}$, where $u_R$ is an affine map taking
the unit cube centered at $0$ to the rectangle $R$; thus $\phi_R$ is
roughly an indicator function of $R$ with ``Schwartz tails".  If ${\cal R}$
is a family of rectangular boxes (usually tubes or plates), we write
$\Phi_{\cal R}=\sum_{R\in{\cal R}}\phi_R$.

We let $\psi(x):\R^{d+1}\rightarrow\R$ be a function with the following properties:

\smallskip

(i) $\psi=\eta^2$, where $\hat{\eta}$ is supported in a small ball centered at $0$.

(ii) $\psi\neq 0$ on a large cube centered at the origin.

(iii) $\sum_{\nu\in\Z^{d+1}}\psi(x-\nu)\equiv 1$.

\smallskip\noindent
We also write $\psi_R=\psi\circ u_R^{-1}$ with $u_R$ as above.  If $R=R_\delta$ is
the unit cube of sidelength $\delta$ centered at $0$, we will write
$\psi_\delta=\psi_{R_\delta}$.

If a family of functions ${\cal F}_\delta$ is given for each $\delta$, we
will say that the functions in ${\cal F}_\d$ are {\em essentially orthogonal} if
$$
\|\sum_{f\in {\cal F}_\d} f\|_2^2\approx \sum_{f\in {\cal F}_\d}\|f\|_2^2.
$$
For example, functions
with finitely overlapping supports or Fourier supports are essentially
orthogonal. 

\medskip

\noindent{\bf Sector decompositions:}

\medskip

$\Pi_{a,\delta}$: $\d$-sectors, defined in (A).

$\M_\delta$: set of centers of $\d$-sectors. For nondegenerate surfaces, we may
take $\M_{\delta}$ to be any $\d^{1/2}$-separated subset of $S$; for conical surfaces, 
$\M_{\delta}$ may be a $\d^{1/2}$-separated subset of $S\cap\{|x|=C_0\}$ for some fixed $C_0$ with
$C_1<C_0<C_2$, where $C_1,C_2$ are as in \bref{def-cone}.
The case of $k$-cones is discussed in Section 7.

$M_\delta=|\M_\delta|$.

$S^\sigma_\d(a)=\sdel\cap\Pi_{a,\sigma}$, for $\d\leq\sigma\lesssim 1$.
Whenever the choice of $a$ is unimportant 
-- which will be most of the time -- we will write
$S^\sigma_\d$ instead of $S^\sigma_\d(a)$.  

$M_{\sigma,\d}$: the number of $\delta$-sectors contained in $S^\sigma_\d$.  

\medskip

Note that 
\begin{equation}\label{eq-M}
M_\d\approx\d^{-(d-k)/2},
\end{equation}
\begin{equation}\label{eq-Msd}
M_{\sigma,\d}\approx M_\delta/M_\sigma\approx(\sigma/\d)^{(d-k)/2}.
\end{equation}
Indeed, the covering and finite overlap conditions in (A) imply that
$M_{\d}|\Pi_{a,\delta}|\approx|\sdel|\approx\d$, hence \bref{eq-M}
follows from \bref{eq-Pi}.  
Also, by the finite overlap and consistency conditions in (A)
we have 
$$|S^\sigma_\d|\approx M_{\sigma,\d}|\Pi_{a,\delta}|
\approx M_{\sigma,\d}|\sdel|/M_\d.$$
In particular, all $S^\sigma_\d$ have approximately the same
size, hence
$$|S^\sigma_\d|\approx |S_\d|/M_\sigma.$$
Comparing the last two estimates and using also \bref{eq-M}, we get \bref{eq-Msd}.

\medskip
\noindent{\bf Dual plates:}
\medskip

We define $\pi^a_0$ to be the rectangular box dual to $\Pi_a-a$:
\[
\pi^a_0=\{x\in\R^{d+1}:\ |x\cdot(y-a)|\leq 1, \ y\in\Pi_a\}.
\]
We fix a tiling of $\R^{d+1}$ by translates of $\pi^a_0$:
$\R^{d+1}=\bigcup_{b\in\Lambda_a}\pi^a_b$, where  $\pi^a_b$
is the translate of $\pi^a_b$ centered at $b$.  We will refer
to $\pi^a_b$, $b\in\Lambda_a$, as {\em plates} dual to the
sector $\Pi_a$.
Thus each $\pi=\pi^a_b$ has a unique sector $\Pi_a$ to which it
is dual; we will sometimes indicate this by writing $\Pi_a=\Pi(\pi^a_b)$.

Note that $\pi^a$ have dimensions $\approx \d^{-1}$, $\d^{-1/2}$,
$1$ in the directions parallel to the short, medium, and long
directions of $\Pi_a$, respectively.


\section{The wave packet decomposition}\label{sec4}


Several basic properties of the norm $\|\cdot\|_{p,\d}$ will be
used throughout this paper. We first record the estimate
\begin{equation}\label{inf.trivial}
\|f\|_{\infty,\d}\lesssim \|f\|_\infty\lesssim M_{\sigma,\d}\|f\|_{\infty,\d},
\ \supp\widehat{f}\subset S^\sigma_\d.
\end{equation}
The first inequality in \bref{inf.trivial}
follows from 
$$\|f\|_{\infty,\d}\leq\max_a\|\Xi_a\ast f\|_\infty
\leq \|\Xi_a\|_1\|f\|_\infty\lesssim \|f\|_{\infty},$$
and the
second one from the identity $f=\sum_{a\in\M_\d} \Xi_a *f$.
Note also that the functions $\Xi_a*f$ have finitely
overlapping Fourier supports and therefore are essentially orthogonal,
so that
\begin{equation}\label{www.98}
\|f\|_{2,\d}^2=\sum_{a\in\M_\d}\|\Xi_a*f\|^2_2 \approx \|f\|_2^2.
\end{equation}

\begin{lemma}\label{interpol}
For all $p\geq 2$ we have
\begin{equation}
\|f\|_{p,\d}\lesssim\|f\|_2^{2/p}\|f\|_{\infty, \d}^{1-2/p}.
\label{looc5}\end{equation}
\end{lemma}

\underline{Proof} Let $f_a=\Xi_a*f$, then
\begin{eqnarray*}
\|f\|_{p,\d}^p &=& \sum_{a}\|f_a\|_p^p
\leq \max_{a}\|f_a\|_\infty^{p-2} \cdot \sum_{a}\|f_a\|_2^2
\\
&=& \|f\|_{\infty,\d}^{p-2}\|f\|_{2,\d}^2.
\end{eqnarray*}
It now suffices to combine this with \bref{www.98}.
\hfill$\square$

We now study the structure of functions with Fourier support in
$\sdel$.  More precisely, we want to decompose such functions into
``Knapp examples", each of which is Fourier localized in a $\d$-sector
and spatially localized (modulo Schwartz tails) in the plate dual to
the sector in question.  The definition below and the next two lemmas
are identical to the corresponding arguments in \cite{Wls}, \cite{LW}
modulo notation and rescaling; we include the proofs for completeness.

\begin{definition}\label{Nfn}
Let $0<\delta\leq\sigma\lesssim 1$.
We define $\Sigma^\sigma_\d(a)$ to be the space of all functions of the form
$f=\sum_{\pi\in\P}f_{\pi},$
where $\P=\P(f)$ is a family of $\d$-plates such that each $\pi^{a'}_b\in\P$
is dual to a sector $\Pi(\pi^{a'}_b)=\Pi_{a'}$ centered at a point ${a'}\in 
S^\sigma_\d(a)$, and
\begin{equation}\label{Nfna}
|f_{\pi}|\lesssim \phi_{\pi},
\end{equation}
\begin{equation}\label{Nfnb}
\supp \widehat{f_{\pi}}\subset \Pi(\pi).
\end{equation}

\end{definition}

If $\P'\subset\P$, we say that $f_{\P'}=\sum_{\pi\in\P'}f_{\pi}$ is 
a {\em subfunction} of $f$.  
When the choice of $a$ is unimportant, we will often omit $a$ from the notation.  Note that
of such functions will be denoted by $\Sigma^\sigma_\d(a)$; as with
for functions in $\Sigma^\sigma_\d$ we have 
$\|f\|_\infty\lesssim M_{\sigma,\d}$ instead of \bref{Nfn-linf}.
If $\sigma\approx 1$, we will sometimes write $\Sigma^\sigma_\d=\Sigma_\d$.

\begin{lemma}\label{Lemma 4.1}
Let $f\in\Sigma^\sigma_\d$.  Then $\supp\widehat{f}\subset S^\sigma_\delta(a)$ and
\begin{equation}\label{Nfn-linf}
\|f\|_\infty\lesssim M_{\sigma,\d},
\end{equation}
\begin{equation}
\|f\|_{p, \d}\lesssim \Big(\sum_{\pi\in\P}|\pi|\Big)^{1/p}
\hbox{ for }2\leq p\leq \infty.
\label{ls3}\end{equation}
\end{lemma}

\underline{Proof} 
The support statement is clear from the definition.  Recall also that
all plates $\pi^a_b $ with fixed $a$ have disjoint interiors; this
immediately implies \bref{Nfn-linf}. It remains to prove \bref{ls3}.
By \bref{looc5}, we only need to do so for $p=2$ and $p=\infty$. 

We first claim that $f_{\pi}$ are essentially orthogonal.  Indeed,
$$
\|f\|_2^2=\Big\|\sum_{a\in\M} \Big(\sum_{b}f_{\pi^a_b}\Big)\Big\|_2^2
\approx 
\sum_{a\in\M} \Big\|\Big(\sum_{b}f_{\pi^a_b}\Big)\Big\|_2^2,
$$
since $\sum_{b}f_{\pi^a_b}$ is Fourier supported in $\Pi_a$, 
and the latter have finite overlap.
It remains to prove the essential orthogonality of $f_{\pi^a_b}$ for
fixed $a$; this follows from an easy argument using the decay
of $\phi_\pi$.
Since $\|f_\pi\|_2^2\lesssim \|\phi_\pi\|_2^2\approx|\pi|$, this yields
\bref{ls3} for $p=2$.

To complete the proof for $p=\infty$ and therefore for all $p$, it
suffices to verify that $\|f* \Xi_a \|_\infty\lesssim 1$
for each $\Xi_a$.  We have 
$f* \Xi_{a}=\sum_{\pi} f_\pi *\Xi_{a}$,
where the only non-zero terms are those corresponding to
$\pi$ whose dual plates $\Pi$ intersect $\Pi_a$.  Since the number of
such $\Pi$ is bounded, and since the plates $\pi$ corresponding to
each $\Pi$ are disjoint, it follows that the plates contributing to
$f*\Xi_a$ have finite overlap.  Hence
\[
\|f* \Xi_a \|_\infty
\lesssim \max_\pi \|f_\pi* \Xi_a \|_\infty\lesssim\max_\pi \|f_\pi\|_\infty
\lesssim 1.
\]
$\ $\hfill$\square$

\bigskip

\begin{lemma}\label{Lemma 4.4}
Assume that $\supp\hat{f}\subset S^\sigma_\d$
and $\|f\|_{\infty,\d}<\infty$. Then there are 
$f_{\l}\in\Sigma^\sigma_\d$, with dyadic $\l$ satisfying
\begin{equation}
\l\lesssim \|f\|_{\infty, \d}
\label{katr1}\end{equation}
and with corresponding plate families $\P_\l$, such that 
\begin{equation}\label{Nfn-dec}
f\approx \sum_{\l}\l f_{\l},
\end{equation}
\begin{equation}
\sum_{\l}\l^p\sum_{\pi\in\P_\l}|\pi|\lesssim\|f\|_{p, \d}^p
\label{wa1}\end{equation}
for each fixed $p\in [2,\infty)$.
\end{lemma}

\underline{Proof} 
We may assume that $\hat f$ is supported in some $\d$-sector $\Pi_a$,
so that $\|f\|_{p,\d}=\|f\|_p$.  Let $\psi_b=\psi_{\pi^a_b}$;
observe that $\|\psi_b f\|_\infty \lesssim \|f\|_{\infty,\d}=\|f\|_\infty$. 
For $\l$ as in \bref{katr1}, we let
$\P_{\l}=\{\pi^a_b:\ \l\leq \|\psi_bf\|_\infty\leq 2\l\}$
and
\[
f_\l=\sum_{\pi_b\in\P_\l}\l^{-1}\psi_{b}^2 f.
\]

To see that $f_{\l}\in\Sigma^\sigma_\d$, it suffices to verify
that each $f_{b}:=\l^{-1}\psi_{b}^2 f$ satisfies \bref{Nfna}, \bref{Nfnb}.
Indeed, \bref{Nfna} is immediate from the definition; also, 
$\widehat{f_b}= \widehat{\psi_b^2}*\widehat f$ and $\widehat{\psi_b^2}$ is
supported in a translate of $c\Pi_a$ centered at $0$ for some $c\ll 1$, hence
\bref{Nfnb} follows.

We have $\sum_\l \l f_\l=\sum_j \psi_b^2 f$ and 
$1\lesssim \sum_b \psi^2_b \lesssim \sum_b \psi_b=1,$
so that \bref{Nfn-dec} follows.  Moreover, by Bernstein's
inequality 
\[
\l^p\approx \|\psi_bf\|^p_{\infty}\lesssim |\Pi_a|\|\psi_bf\|_p^p
\approx |\pi^a_b|^{-1}\|\psi_bf\|_p^p,
\]
hence
\[
\sum_{\l}\l^p\sum_{b:\pi^a_b\in\P_\l}|\pi_b|
\lesssim \sum_b\|\psi_bf\|_p^p\lesssim\|f\|_p^p
\]
as required.\hfill$\square$


\section{Proof of Theorem \ref{thm1}}\label{sec5}


Let $p$ be as in Theorem \ref{thm1}. 
In the inductive argument that follows, we will need a
variant of our main estimate for functions with Fourier support
in sectors $\sdel^\sigma$.  In this setting, our main inequality may be stated as
follows: suppose that $\supp\widehat{f}\subset S_\delta^\sigma$ and
\begin{equation}\label{www.99}
\|f\|_{p, \d}\leq 1,
\end{equation}
then
\begin{equation}\label{www.100}
\|f\|_p^p \lesssim \d^{-\epsilon}M_{\sigma,\delta}^{rp+\epsilon}
\end{equation}
for any $\epsilon>0$, where we recall that 
\begin{equation}\label{def-r}
r=1-\frac{2}{p}-\frac{2}{p(d-k)}.
\end{equation}

\begin{proposition}\label{ttt.prop1}
For any $\e>0$, we have the estimate
\begin{equation}
\|f\|_p^p \lesssim \d^{-\epsilon}M_{\sigma,\delta}^{rp+\epsilon}\|f\|_2^2,
\label{mnmn}\end{equation}
for all $f$ with $\|f\|_{\infty,\d}\lesssim 1$ and $\supp\widehat f\subset
S^\sigma_\delta$.
\end{proposition}

Assuming Proposition \ref{ttt.prop1}, let us complete the proof of
\bref{www.100}.  Pich $f$ with $\supp\widehat{f}\subset S_\delta^\sigma$ and
obeying \bref{www.99}, then by Lemma \ref{Lemma 4.4} we have
\[
f\approx \sum_{\l\leq \|f\|_\infty}\l f_{\l},\ f_{\l}\in\ssd.
\]
From \bref{wa1} and Lemma \ref{Lemma 4.1} we have
$\|f_{\l}\|_{\infty,\d}\lesssim 1$ and 
\begin{equation}\label{www.97}
1\gtrsim \l^p\sum_{\pi\in \P_\l}|\pi|\gtrsim
\l^p\|f_\l\|_{2,\d}^2\approx \l^p \|f_{\l}\|_2^2,
\end{equation}
which together with \bref{mnmn} yields \bref{www.100} for $\l f_\l$.

The conclusion follows by summing over $\l$ if we can show that the
summation can be restricted to the logarithmically many values of $\l$
in $[\d^{K},\d^{-K'}]$ for some $K,K'$.
Indeed, for any $C$ we have
$\sum_{\l\leq\d^{K}}\l f_{\l}\lesssim \d^C$ if $K$ is large
enough, hence \bref{www.100} holds for this part.
It remains to prove that $\|f\|_\infty\lesssim \d^{-K'}$.
Let $f_a=\Xi_a *f$, then
\[
\|f_a\|_\infty\lesssim |\Pi_a|^{1/p}\|f_a\|_p
\lesssim |\Pi_a|^{1/p},
\]
by Bernstein's inequality and \bref{www.99}.
The desired bound follows on summing over $a$.
\qed

\bigskip

Proposition \ref{ttt.prop1} is an immediate consequence of Proposition
\ref{cor-iterate} and Corollary \ref{strong-type} below.  
The main inductive argument is given in Proposition \ref{scales}, the
proof of which will occupy Sections 4-6.

\begin{definition}\label{def-P}
We say that $P(p,\alpha,\e)$ holds if for all functions $f$ such that
$\supp\widehat{f}\subset S^\sigma_\d$ and $\|f\|_{\infty, \d}\leq 1$ we have
\begin{equation}
|\{ |f|>\l\}|\lesssim \l^{-p} \d^{-\epsilon}M_{\sigma,\delta}^{rp+\alpha}\|f\|_2^2\ ,
\label{4.1}\end{equation}
for all $0<\d\leq\sigma\lesssim 1$, provided that $\d$ is small enough. 
\end{definition}

\begin{proposition}\label{cor-iterate}
$P(p,\alpha,\epsilon)$ holds for all $\alpha,\e>0$.
\end{proposition}

\underline{Proof} We will assume that $f$ is as in Definition \ref{def-P}.   
Observe that \bref{4.1} follows 
automatically from Chebyshev's inequality if
\begin{equation}\label{cheb}
\l^{p-2}\lesssim \d^{-\e}M_{\sigma,\delta}^{rp+\alpha},
\end{equation}
By \bref{inf.trivial}, we may assume that
$\l\leq \|f\|_\infty\lesssim M_{\sigma,\delta}$.
Thus for \bref{cheb} to hold, it suffices if
$M_{\sigma,\d}^{p-2}\lesssim \d^{-\e}M_{\sigma,\delta}^{rp+\alpha}$,
or in other words
$$
M_{\sigma,\d}^{p-2-rp-\alpha}=M_{\sigma,\delta}^{\frac{2}{d-k}-\alpha}
\lesssim\d^{-\e}.
$$
This has two consequences of interest to us:

\smallskip

1. \bref{4.1} holds for all $\alpha >0$ if $\d^{-\e}
\gtrsim M_{\sigma,\delta}^{\frac{2}{d-k}}\approx \sigma/\delta$  (the last
equality comes from \bref{eq-Msd}),

\smallskip
2. $P(p,\alpha,\e)$ holds for any $\e>0$ if $\alpha\geq\frac{2}{d-k}$, since
then $M_{\sigma,\delta}^{\frac{2}{d-k}-\alpha}\lesssim 1$.

\smallskip

The main inductive step is contained in the following proposition.

\begin{proposition}\label{scales}
Fix $p>p_d$ and $0<\e_2<\e$.  Suppose that $P(p,\alpha,\e)$ holds,
and let $\e_0>0$ be sufficiently small (depending only on $p,d,k$).
Then \bref{4.1} holds with $\alpha$ and $\e$ replaced by 
$(1-\frac{\e_0}{6})\alpha$ and $4\e$, respectively, provided
that $\sigma/\d\gtrsim \d^{-\e_2}$.
\end{proposition}

Assuming Proposition \ref{scales} for the moment, let us finish
the proof of Proposition \ref{cor-iterate}.
Fix $\alpha,\e>0$.  Fix also $\alpha_0>2/(d-k)$, 
sufficiently small positive numbers $\e_1$, $\e_2$ 
(depending on $d,k,p$), and a large integer $m_0$ so that
\begin{equation}\label{epses}
m_0>\frac{\log(\alpha/\alpha_0)}{\log(1-\e_0/6)},\ 
\e_1<4^{-m_0}\e,\ \e_2<\e_1.
\end{equation}
By 2. above, $P(p,\alpha_0,\e_1)$ holds.  

Assume that we know
$P(p,(1-\e_0/6)^m\alpha_0,4^m\e_1)$ for some integer $m\geq 0$.  We 
claim that this implies $P(p,(1-\e_0/6)^{m+1}\alpha_0,4^{m+1}\e_1)$.
Indeed, if 
$$\sigma/\d\lesssim  \d^{-4^{m+1}\e_1},$$
then \bref{4.1} follows from 1. above.  Otherwise, we must have
$$\sigma/\d\gtrsim  \d^{-4^{m+1}\e_1}
\geq  \d^{-\e_1}\geq\d^{-\e_2},$$
in which case the claim follows from Proposition \ref{scales}.
After $m_0$ iterations, we obtain
$P(p,(1-\e_0/6)^{m_0}\alpha_0,4^{m_0}\e_1)$.  From
\bref{epses} we have
$(1-\e_0/6)^{m_0}\alpha_0<\alpha$ and $4^{m_0}\e_1<\e$, hence
$P(p,\alpha,\e)$ follows as required.
\qed

\begin{corollary}\label{strong-type}
Assume that $P(p,\alpha,\e)$ holds.  Then for all $f$ with
$\supp\hat{f}\subset S^\sigma_\d$ and $\|f\|_{\infty,\d}\lesssim 1$ we have
\begin{equation}
\|f\|_p^p \lessapprox \d^{-\epsilon} M_{\sigma,\delta}^{rp+\alpha}\|f\|_2^2.
\label{mnmnmn}\end{equation}
\end{corollary}

\underline{Proof}
Write $|f|\approx \sum_\l \l\chi_{\{|f|\approx\l\}}$ with dyadic $\l$.
Then $\|f\|_\infty\lesssim M_{\sigma,\delta}$ by \bref{inf.trivial};
also, \bref{mnmnmn} is trivial if 
$\|f\|_\infty\lesssim\d^{n}$ for $n$ large enough.
The lemma follows by summing \bref{4.1} over dyadic $\l$ 
with $\d^n\lesssim \l\lesssim M_{\sigma,\delta}$.
\qed

\bigskip


\section{A substitute for scaling}\label{sec-scale}


In this section we develop the geometrical arguments which will replace the 
scaling arguments of \cite{Wls},
\cite{LW}.  Instead of rescaling $S$ by powers of 
$\delta$, we keep $S$ fixed and consider its $\rho$-neighbourhoods for
intermediate values of $\rho$ between $\delta$ and 1.  We then need a 
mechanism for efficient conversion of functions with Fourier support in $\sdel$
to functions with Fourier support roughly equal to $S_\rho$.  This is done as follows.
Take a function $f$ with $\supp\widehat f\subset \sdel$, and let
$f_Q=\psi_Q f$, where $Q$ is a $\rho^{-1}$-cube. 
This localizes $f$ spatially in $Q$, modulo Schwartz tails, and (since
$\widehat{f_Q}=\widehat{\psi_Q}*\widehat{f}$) extends its Fourier support 
to $S_\rho$.  It is instructive to examine what happens to a function
$f=f_\pi$ satisfying \bref{Nfna}, \bref{Nfnb}.  
If $\rho\lesssim\sqrt{\d}$ (which we will assume in all applications of
the lemma), $f_Q$ is essentially obtained from $f$ by shortening the
supporting plate $\pi$ to length about $\rho^{-1}$ in the long
direction; to compensate for it,
the thickness of the Fourier support $\Pi$ increases to $\rho$.

\begin{lemma}\label{Lemma 4.2}
Assume that $\supp \widehat{f}\subset \sdel^\sigma$, and let $Q$ be a $\rho^{-1}$-cube
for some $\d\lesssim\rho\lesssim 1$. Let $f_Q=\psi_Qf$.  Then
$\supp\widehat{f_Q}\subset S_{\rho}^\sigma$ and
\begin{equation}
\|f_Q\|_{p,\rho}\lesssim (\d/\rho)^{1/p}M_{\rho,\d}^{1-2/p}
\|f\|_2^{2/p}\|f\|_{\infty, \d}^{1-2/p}
\label{ls20}\end{equation}
for all $p\geq 2$.
\end{lemma}

\underline{Proof}
Observe that $\widehat{\psi_Q}$ is supported in a $\rho$-cube
centered at 0, and has size about $\rho^{-d-1}$ on a proportional cube.
We have $\widehat{f_Q}=\widehat{\psi_Q}*\widehat{f}$.  
The support statement follows immediately as discussed above.
By \bref{looc5}, it suffices to prove \bref{ls20} for
$p=2$ and $p=\infty$.  We have $f_Q=\psi_Qf=\sum_a\psi_Q(\Xi_a*f)$.
If we convolute this with a function whose Fourier transform 
is supported in a $\rho$-sector, the only contributing 
$\delta$-sectors will be those that intersect the $\rho$-neighbourhood 
of the $\rho$-sector in question, hence the $L^\infty$ bound follows.
For $p=2$, we write 
\[
\widehat{f_Q}(\xi)=\widehat{\psi_Q}*\widehat{f}(\xi)
=\int\widehat{\psi_Q}(\xi-\eta)\chi_{\sdel}(\eta)
\widehat{f}(\eta)d\eta,
\]
and use Schur's test.  We have
$\int\widehat{\psi_Q}(\xi-\eta)\chi_{\sdel}(\eta)
d\eta\lesssim \d/\rho$
and
$\int\widehat{\psi_Q}(\xi-\eta)\chi_{\sdel}(\eta)
d\xi\lesssim 1.$  Hence
\begin{equation}\label{Schur}
\|f_Q\|_2 =\|\widehat{f_Q}\|_2
\lesssim (\d/\rho)^{1/2}\|\widehat{f}\|_2
= (\d/\rho)^{1/2}\|{f}\|_2
\end{equation}
as required.  
\hfill$\square$

\bigskip

We also need a substitute for the Lorentz transformations used
in \cite{Wls}, \cite{LW}, \cite{PS}.  The key geometrical observation
turns out to be the following.  Suppose that $\widehat f$ is supported
only in a small part of $\sdel$, say in $S^\sigma_\delta$ for 
some $\delta \lesssim \sigma\ll 1.$  Compared to $\sdel$, $S^\sigma_\delta$ 
is quite flat.  We may therefore convolute $\widehat{f}$ with the
characteristic function of a box ${\cal R}$
whose dimensions in these ``flat" directions are larger than $\rho$,
and still stay in a $\rho$-neighbourhood of $S^\sigma_\d$.  For 
$\sigma\ll 1$, this box can be quite a bit larger than a $\rho$-cube;
in the special case when $\sigma=\rho=\delta$, it will be a
$\d$-sector as opposed to a $\d$-cube.  This will result in a
considerable gain in \bref{ls30} below.  The reader may be interested to
verify that replacing $|R|$ by $|Q|=\rho^{-d-1}$ in \bref{ls30}
would have disastrous consequences at the end of the proof of
Theorem \ref{thm1}.

\begin{lemma}\label{box-convolute}
Let $\delta\lesssim\rho\lesssim \sigma \lesssim 1$, and assume
that $\sigma\leq \d^{1/2}$.  
Let $f$ satisfy $\supp \widehat{f}\subset S^\sigma_\delta(a)$. 
We define $\calr$ to be a box centered
at 0, with sides parallel to the sides of $\Pi_a$ and sidelengths
$\rho,\rho\sigma^{-1/2},\rho\sigma^{-1}$ in directions parallel
to the short, medium, and long directions of $\Pi_a$ respectively.
We also let $R_0$ be the box dual to $\calr$ centered at 0.  Fix
a translate $R$ of $R_0$, and let $f_R=\psi_R f$.  Then
\begin{equation}\label{R-eq1}
|R|=|\calr|^{-1}=\rho^{-d-1}\sigma^{(d+k)/2},
\end{equation}
and:

\smallskip
(i) $\supp\widehat{f_R}\subset S^\sigma_\rho$,

\smallskip
(ii) $f_R$ obeys \bref{ls20} with $p=\infty$, i.e. 
$\|f_R\|_{\infty,\rho}\lesssim M_{\rho,\d} \|f\|_{\infty, \d}$,

\smallskip
(iii) if $\rho\lesssim\sqrt{\sigma\d}$, we have the estimate
\begin{equation}
\|f_R\|_2\lesssim M_{\sigma,\delta}^{1/2}|R|^{1/2}\|f\|_{\infty, \d}.
\label{ls30}\end{equation}
\end{lemma}

Note that if $\rho\approx\sqrt{\sigma\d}$, then $R$ has dimensions
$(\sigma\d)^{-1/2},\d^{-1/2},(\sigma/\d)^{1/2}$, and
\bref{R-eq1} becomes
\begin{equation}\label{R-eq2}
|R|\approx \rho^{(k-1)/2}\d^{-(d+1)/2}.
\end{equation}

\medskip

\underline{Proof} We have $\widehat{f_R}=\widehat{\psi_R}*\widehat{f}$
and $\supp\widehat{\psi_R}\subset {\calr}$.
Thus the first two parts of the lemma follow by the same
argument as in the proof of Lemma \ref{Lemma 4.2}, if we can show that
for all $b\in S\cap S_\d^\sigma(a)$:
\begin{equation}\label{box-e1}
\calr+b\subset S_\rho^\sigma(a)
\end{equation}
and
\begin{equation}\label{box-e2}
S_\d^\rho(b)+\calr\subset \Pi_{\rho,a}.
\end{equation}
Indeed, for any such $b$ we have
$$
\Pi_{\sigma,b}\subset C\Pi_{\sigma,a}\subset C'\Pi_{\sigma,b}.
$$
It follows that 
$\Pi_{\sigma,a}+(b-a)\subset C\Pi_{\sigma,b}$, hence
$\rho\sigma^{-1} \Pi_{\sigma,a}+(b-a)\subset C\rho\sigma^{-1} \Pi_{\sigma,b}$.  
But $\rho\sigma^{-1} \Pi_{\sigma,a}-a=\calr$, and 
$C\rho\sigma^{-1} \Pi_{\sigma,b}\subset C'\Pi_{\rho,b}$, by comparing
sidelengths and using that $\rho\lesssim\sigma$.  Thus 
\begin{equation}\label{box-e3}
\calr+b\subset C'\Pi_{\rho,b},
\end{equation}
which immediately implies both \bref{box-e1} and \bref{box-e2}.

It remains to prove (iii).  
If $\rho\lesssim \sqrt{\sigma\d}$, the same argument as in the proof
of \bref{box-e3} shows that $b+\calr$ is contained in the box ${\cal B}_b$ obtained
from $C\Pi_{\delta,b}$ by thickening its shortest sidelength to $C\rho$.
Hence the functions on the right-hand side of the identity
$f_R=\sum_{a\in\M} \psi_R \cdot(\Xi_{a}\ast f)$ are essentially orthogonal
since their Fourier supports have finite overlap.   It follows that
\[
\|f_R\|_2^2\lesssim\sum_{a}\|\psi_R\cdot(\Xi_a\ast f)\|_2^2
\lesssim\sum_{a}\|f\|_{\infty, \d}^2\,\|\psi_R\|_2^2
\lesssim M_{\sigma,\d}|R|\|f\|_{\infty, \d}^2
\]
as claimed.
\hfill$\square$


\section{The localization property}\label{sec2}


Throughout this section we will assume that $\sigma/\d\gtrsim\d^{-\e_2}$,
where $\e_2>0$ was fixed in Section \ref{sec5}.

\begin{definition}\label{localize}
Let $f\in\ssd$. We say that $f$ {\em localizes} at $\lambda$
if there are subfunctions $f^Q$ of $f$, where $Q$ runs over $t\d^{-1}$-cubes
(recall that $t\approx(\sigma/\d)^{\e_0}$), 
such that
\begin{equation}
\sum_Q|\P(f^Q)|\lessapprox |\P(f)|
\label{vo9}\end{equation}
and
\begin{equation}\label{zzz.45} 
|\{|f|\geq\l\}|\lessapprox  \sum_Q |Q\cap \{|f^Q|\gtrapprox \l\}|.
\end{equation}
\end{definition}

Our task is now to find conditions under which $f$ localizes.

\begin{lemma}\label{Lemma 2.3}
Let $\P$ be a family of plates, and let $W\subset \R^{d+1}$.  
Then there is a relation
$\sim$ between plates in $\P$ and $t\d^{-1}$-cubes $Q$ such that
\begin{equation}\label{propR}
|\{Q:\ \pi\sim Q\}|\lessapprox 1\hbox{ for all }\pi\in\P,
\end{equation}
and
\begin{equation}\label{badI}
I_b\lessapprox t^{-c_1}|W||\P|^{\frac{1}{2}},
\end{equation}
where
\[
I_b=\int_W \sum_{\pi\in\P,\pi\not\sim Q(x)}\chi_\pi(x)
=\sum_{\pi\in\P}|\{x\in W\cap\pi:\ Q(x)\not\sim \pi\}|.
\]
\end{lemma}

\underline{Proof} We may assume that $W\subset\{|x|\lesssim \d^{-1}\}$.
For each $\pi\in\P$, we let $Q(\pi)$ be the $t\d^{-1}$-cube with maximal
$|W\cap Q\cap\pi|$ (if there is more than one, pick one arbitrarily). 
We then say that $\pi\sim Q$ if $Q\cap 10 Q(\pi)\neq\emptyset.$  
Clearly, the number of such cubes (for a fixed $\pi$) is $\lesssim 1$.
We now prove \bref{badI}.

By dyadic pigeonholing, there are $\nu$ and $\P'\subset \P$ such that
$|I_b|\lessapprox \nu |\P'|$ and 
\[
|\{x\in W\cap\pi:\ \pi\not\sim Q(x)\}|\approx\nu\hbox{ for each }\pi\in \P'.
\]
Thus for each $\pi\in\P'$ there is a cube $Q'(\pi)$ such that $\pi\not\sim
Q'(\pi)$ and $|W\cap Q'(\pi) \cap \pi|\gtrsim t\nu$.  But then 
$|W\cap Q(\pi)\cap \pi|\gtrsim t\nu$, by the definition of $Q(\pi)$.
The total number of cubes covering $W$ is $\lesssim t^{-d-1}$, hence we may choose $Q,Q'$
so that $Q=Q(\pi)$ and $Q'=Q'(\pi)$ for at least $t^{2d+2}|\P'|$ plates $\pi\in\P'$.
Let
\[
A=\sum_{\pi\in\P'}|W\cap Q\cap \pi|\cdot |W\cap Q'\cap\pi|,
\]
then for $Q,Q'$ as above we have
\[
A\gtrsim t^{2d+2}|\P'|\cdot (t\nu)^2=t^{2d+4}\nu^2|\P'|.
\]
But on the other hand,
\[
A=\int_{W\cap Q}\int_{W\cap Q'}\sum_{\pi\in\P'}\chi_\pi(x)\chi_\pi(x')dx'dx\]
\[
=\int_{W\cap Q}\int_{W\cap Q'}|\{\pi\in\P':\ x,x'\in\pi\}|dxdx'.
\]
We claim that the integrand is bounded by $t^{-d}$. Indeed, let
$x\in Q$, $x'\in Q'$, then $|x-x'|\gtrsim t\d^{-1}$.  If $x,x'\in\pi^a_b$. 
then the angle between $x-x'$ and $n(a)$ is $\lesssim t^{-1}\d^{1/2}$;
but by (A), $|n(a)-n(a')|\geq\d^{1/2}$ for $a\neq a'$, hence there are
at most $\lesssim t^{-d}$ distinct $a$'s as above.  Thus
\[
A\lesssim t^{-d}|W\cap Q|\cdot|W\cap Q'|
\lesssim t^{-d}|W|^2.
\]
From this and the lower bound just stated, we have\footnote{
We could keep track of the exact powers of $t$ but will have no need
to do so.}
$\nu\lesssim t^{-c}|W|\,|\P'|^{-1/2}$, so that
\[
|I_b|\lessapprox \nu|\P'|\lesssim t^{-c}|W|\,|\P'|^{1/2}
\leq t^{-c}|W|\,|\P|^{1/2}
\]
as required.
\hfill$\square$

~

In the sequel, we will use a version of Lemma \ref{Lemma 2.3} with 
``Schwartz tails".  The proof follows word-for-word the standard
argument given in detail in \cite{Wls} and \cite{LW}, and we see no 
reason to reproduce it here.

\begin{lemma}\label{Lemma 2.7} Let $W\subset\R^{d+1}$ be measurable, 
and let $\P$ be a family of plates. Fix a large constant $M_0$. Then, if
the constant $K$ in the definition of $\phi$ has been chosen large
enough, there is a relation $\sim$ between $t\d^{-1}$-cubes $Q$ and 
plates in $\P$ satisfying \bref{propR} and such that
\begin{equation}
\int_W\Phi_{\P}^b\lessapprox t^{-c}|\P|^{\frac{1}{2}}|W| +\d^{M_0}|W|,
\label{vu2}\end{equation}
where
\begin{equation}
\Phi^b_{\P}(x)=\sum_{\pi\in\P,\pi\not\sim Q(x)}\phi_{\pi}(x)
\label{sch2}\end{equation} 
\end{lemma}

\begin{lemma}\label{Lemma 3.1}
Let $f\in\ssd$, with plate family $\P$. Assume that 
\begin{equation}
|\P|\leq t^{2c+2}\l^2
\label{vo8}\end{equation} 
with $c$ as in \bref{vu2}. Then $f$ localizes at $\l$.
\end{lemma}

\underline{Proof} Let $W=\{|f|\geq\l\}$,
and let $\sim$ be the relation defined in \ref{Lemma 2.7}.  
For each $Q$, let $f^Q=\sum_{\pi\sim Q}f_\pi$. By \bref{propR},
\bref{vo9} holds.  Also, we have
\[
\int_{W}\Phi_{\P}^b\lessapprox t^{-c}|\P|^{\frac{1}{2}}|W|
\lesssim t\l|W|.
\]
Hence we must have $\Phi_{\P}^b(x)\lesssim t\l$ on some set
$W^*\subset W$ with proportional measure.
Let $x\in W^*\cap Q$, then
\[
|f(x)-f^Q(x)|=|\sum_{\pi\not\sim Q}f_\pi(x)|
\lesssim \Phi_{\P}^b(x)\lesssim t\l,
\]
so that $|f^Q(x)|\gtrsim\l$ as claimed.
\hfill$\square$

~

\begin{lemma}\label{Lemma 3.2} 
Let $f\in\ssd(a)$, with plate family $\P$.  Fix
a tiling $\{R\}$ of $\R^{d+1}$ by translates of the 
rectangular box $R_0$ defined in Lemma \ref{box-convolute}, with
$\rho\approx \sqrt{\sigma\delta}$. 
Then either $f$ localizes at $\l$, or else there is a subfunction
$f^*$ of $f$ such that $|f^*|\gtrapprox \l$ on a logarithmic fraction
of $\{|f|\geq\l\}$ and 
\begin{equation}
\|\psi_{R}f^*\|_{2}^2
\lesssim t^{-C}\l^{-2}(\sigma/\d)^{-k/2}M_{\sigma,\d}|R||\P|.
\label{vo5.1}\end{equation}
for each $R$.

\end{lemma}

\underline{Proof}
For each $\pi=\pi^b_{b'}\in\P$ we let $\tau^b_{b'}$ be a rectangular
box containing $\pi^b_{b'}$ with the same dimensions except that all short
sidelengths are extended from 1 to $\sqrt{\sigma/\delta}$.  From this
family of boxes, we choose a maximal subset $\tilde\T$ with the property
that $\tau'\not\subset C\tau$ for any $\tau',\tau\in\tilde\T$ and for
a suitable large constant $C$, and let $\T=\{2C\tau:\tau\in\tilde\T\}$.  
Abusing the notation and standard terminology, we will refer to the boxes
in $\T$ as {\em tubes} and continue to denote them by $\tau$.
Each $\pi\in\P$ is then contained in some $\tau(\pi)\in\T$ (if there is
more than such $\tau$, we choose one arbitrarily).
Observe that the largest angle between two line segments of length
$\sim\d^{-1}$ contained in a tube $\tau$ is bounded by
$\d\sqrt{\sigma/\d}\lesssim\sqrt{\d}$, hence the angle separation
condition in (A) implies that 
\begin{equation}\label{evr1}
|\{b\in\M_\d:\ \pi_{b'}^b\subset\tau\hbox{ for some }b'\}|\lesssim 1.
\end{equation}

We let $\T_m$ be the set of those $\tau\in\T$ with
$|\{\pi\in\P:\ \tau(\pi)=\tau\}|$ between $m$ and $2m$.
Let also $\P_m=\{\pi\in\P:\ \tau(\pi)\in\P_m\}$.
Then for some $m$ we have 
\[
|W|\gtrapprox |\{|f|\geq\l\}|,
\]
where $W=\{x:\ |f(x)|\geq\l,\ |f_{\P_m}(x)|\gtrapprox \l\}$.
Fix this value of $m$, then
\[
|\T_m|\lessapprox |\P|/m.
\]
We now consider two cases. 

~

{\it Case 1:} If $\l\geq t^{-C}(|\P_m|/m)^{1/2}$,
$f$ localizes. The proof is similar to that of
Lemma \ref{Lemma 3.1}: we first construct a relation $\sim$
between $t\d^{-1}$-cubes and tubes in $\T_m$ such that \bref{propR} holds, and
\begin{equation}
\int_W\Phi_{\T_m}^b\lessapprox t^{-C}|\T_m|^{1/2}|W| +\d^{M_0}|W|,
\label{vu2a}\end{equation}
where $M_0$ is a large constant and
\begin{equation}
\Phi^b_{\T_m}(x)=\sum_{\tau\in\T_m,\tau\not\sim Q(x)}\phi_{\tau}(x).
\label{sch2a}\end{equation} 
The construction is identical to that in Lemma \ref{Lemma 2.7}, therefore
we omit it.  Since $|\T_m|\lessapprox |\P_m|m^{-1}$, it follows that
\[
\int_{W}\Phi_{\T_m}^b
\lessapprox t^{-C} \Big(\frac{|\P_m|}{m}\Big)^{1/2}|W|
\lesssim t\l |W|.
\]
Hence $\Phi_{\T_m}^b\lesssim t\l$ on a subset $W^*\subset W$ with proportional measure.
We write $\pi\sim Q$ if $\tau(\pi)\sim Q$, and let
\[
f^Q=\sum_{\pi\sim Q}f_{\pi}.
\]
Then \bref{vo9} follows from \bref{propR}, and
for $x\in W^*\cap Q$ we have
\[
|f(x)-f^Q(x)|=|\sum_{\pi\not\sim Q}f_\pi(x)|
\lesssim \Phi_{\T_m}^b(x)\lesssim t\l,
\]
so that $|f^Q|\gtrsim\l$ on $W^*\cap Q$ as required.
\hfill$\square$

~ 

{\it Case 2:} Assume now that $\l\leq t^{-C}(|\P_m|/m)^{1/2}$; we
will show that $f_{\P_m}$ satisfies \bref{vo5.1}.  
Fix $R$. 
It is easy to see (cf. the proof of Lemma \ref{Lemma 4.1})
that $\psi_{R}f_{\pi}$ are essentially orthogonal, hence
\[
\|\psi_{R}f_{\P_m}\|_2^2\approx 
\|\sum_{\pi\in\P_m}\psi_{R}f_{\pi}\|_2^2
\lesssim\int\sum_{\pi\in\P_m}|f_{\pi}|^2\phi_{R}
\lesssim\int\Phi_{\P_m}\phi_{R},
\]
where at the last step we used \bref{evr1}.

Recall that if $\rho=\sqrt{\sigma/\d}$ then $R$ has dimensions $(\sigma
\d)^{-1/2}$, $\d^{-1/2}$, $(\sigma/\d)^{1/2}$.  
The main geometrical observation is that if a plate $\pi=\pi^b_{b'}$
intersects $CR$ at all, then a piece of $\tau(\pi)$
of length $\sim (\sigma\delta)^{-1/2}$ is entirely
contained in $C'R$ for a suitable choice of $C'$.  Indeed, let ${\cal B}_b$ be 
the box defined in the proof of Lemma \ref{box-convolute}, and recall the
inclusion $\calr+b\subset {\cal B}_b$.  Thus the converse inclusion holds
for the dual boxes.  Hence if $\pi$ intersects $CR$, 
then a piece of $\pi$ of length $\sim (\sigma\delta)^{-1/2}$ (which is
dual to ${\cal B}_b$) is entirely contained in $C'R$.
Since all dimensions of $R$ are $\gtrsim (\sigma/\d)^{1/2}$, we may
increase $C'$ (if necessary) to obtain the same inclusion for $\tau$.

The corresponding Schwartz tails estimate is
\begin{equation}\label{xxx01}
\int_R \sum_{\pi:\tau(\pi)=\tau} \phi_{\pi}\lesssim 
m(\sigma/\delta)^{-k/2}\int_{C'R}\phi_{\tau}
\end{equation}
for all $R$ and $\tau\in\T$.  Namely, 
let $T$ be the infinite tube extending $\tau$ in the direction of its
longest axis.  If $R\cap CT\neq\emptyset$, \bref{xxx01} follows from the
above observation and from the fact that $|\tau|=(\sigma/\d)^{k/2}|\pi|$.
Otherwise, we have the pointwise bound
$\sum_{\pi:\tau(\pi)=\tau} \phi_{\pi}(x)\lesssim 
m(\sigma/\delta)^{-k/2}\phi_{\tau}(x)$, which again yields \bref{xxx01}.

Summing \bref{xxx01} over all $\tau\in\T$, we obtain 
\[
\int_R\Phi_{\P_m}\lesssim m(\sigma/\delta)^{-k/2}\int_{C'R}\Phi_{\T}
\]
for all $R$.  By an easy covering argument,
\[
\int\Phi_{\P_m}\phi_{R}
\lesssim (\sigma/\d)^{-k/2} m\int\Phi_{\T}\phi_{C'R}.  
\]
Fix a point $x$, then the number of tubes with $x\in\tau$ is bounded by
$CM_{\sigma,\d}$.  We convert this to a Schwartz tails bound
$\Phi_{\T_m}\lesssim M_{\sigma,\d}$, and deduce that 
\[
\int\Phi_{\P_m}\phi_{R}
\lesssim m (\sigma/\d)^{-k/2}M_{\sigma,\d}\int\phi_{C'R}
\lesssim m (\sigma/\d)^{-k/2}M_{\sigma,\d}|R|.
\]
Recall also that $m\leq t^{-C} |\P_m|\l^{-2}$.  Therefore
\[
\|\psi_{\Delta}f_{\P_m}\|_2^2
\lesssim\int\Phi_{\P_m}\phi_{R}
\lesssim t^{-C}|\P|\l^{-2}(\sigma/\d)^{-k/2}M_{\sigma,\d}|R|.
\]
\hfill$\square$


\section{Proof of Proposition \ref{scales}}\label{sec5a}


In this section we complete the inductive argument.  The outline of the
proof is as follows.  Let $\supp\widehat f \subset S^\sigma_\delta$;
we want to estimate $\|f\|_p$ in terms of $\|f\|_{p,\delta}$.
To this end, we introduce an
intermediate scale $\rho=\sqrt{\sigma\d}$ and first break up
$S^\sigma_\delta$ into smaller pieces $S^\rho_\delta$, then
decompose these further into $\d$-sectors.  On each $S^\rho_\delta$
we apply the inductive assumption $P(p,\alpha,\e)$.  We would now
like to complete the proof by applying $P(p,\alpha,\e)$ on $S^\sigma_\rho$
to deal with the coarse decomposition.  This, however, would not
increase the value of $\alpha$; therefore at this point we want
to change scales from $\rho$ to $\rho^{1-\e_0}$, which will ensure
the desired gain.  We will see that this is possible if $f$ localizes.
Hence the coarse-scale decomposition in Lemma \ref{Lemma 5.1} is designed
so as to allow localization, either on scale $\rho$ or on a second
intermediate scale between $\rho$ and $\sigma$.

We continue to assume that 
\begin{equation}\label{sd-size}
\sigma/\delta\gtrsim\d^{-\e_2},
\end{equation}
where we recall that $\e_2$ was chosen in Section \ref{sec5} so that 
$\e_2<\e$.  Fix also a small positive
number $\e_3$ with $\e_3<e_2^2$.

\begin{lemma}\label{Lemma 5.1}  Assume that $P(p,\alpha,\e)$ is known
for some $p$ and $\alpha$.
Assume that $\supp\widehat{f}\subset S^\sigma_\d$ and
that $\|f\|_{\infty, \d}\lesssim 1$. Let $\rho=\sqrt{\sigma\d}$.  Let
also $\calr, R_0$ be the boxes defined in Lemma \ref{box-convolute}, and fix
a tiling $\{R\}$ of $\R^{d+1}$ by translates of $R_0$.
Then for any $\l\geq\d^{C}$ and for any $\e>0$ we may find a $\l_*$ 
and functions $f_{R}\in\Sigma^\sigma_\rho$, 
with respective plate families $\P_R$, such that
a logarithmic fraction of $\{|f|\geq\l\}$ is contained in
$\bigcup_R \{|f_R|\geq\l_*\}$ and
\begin{equation}\label{vp.lambda}
\d^{2\epsilon_3} \l M_{\rho,\d}^{-1}\lesssim \l_*\lesssim M_{\sigma,\rho},
\end{equation}
\begin{equation}
\sum_{\pi\in\P_{R}}|\pi|\leq\d^{-C\e_3}(\frac{\l_*}{\l})^2
\|\psi_R f\|_2^2,
\label{vp2}\end{equation}
\begin{equation}
\sum_{R}\sum_{\pi\in\P_{R}}|\pi|\leq
(\frac{\l_*}{\l})^{p}
\d^{-\e-C\e_3} M_{\rho,\d}^{rp+\alpha}\|f\|_2^2.
\label{vp3}\end{equation}
\end{lemma}

\underline{Proof} We write $f=\sum_R \psi_Rf$. 
Using \bref{inf.trivial} and the Schwartz decay of $\psi$,
it is easy to prove that
\begin{equation}
\{|f|\geq\l\}\subset\bigcup_R \{|\psi_R f|\geq c \d^{\e_3}\l\}.
\label{ls13}\end{equation}

By Lemma \ref{box-convolute}, $\supp\widehat{\psi_Rf}\subset S^\sigma_{\rho}$ and 
$\|\psi_Rf\|_{\infty,\rho}\lesssim M_{\rho,\d}$.
Lemma \ref{Lemma 4.4} now yields a decomposition
\[
\psi_R f \approx \sum_{h}hg^R_h,
\]
where $g^R_h\in\Sigma^\sigma_\rho$ and $h\lesssim M_{\rho,\d}$.

Since we assume that $\l\geq\d^{C}$, there are logarithmically many 
relevant dyadic values of $h$.  We may therefore choose $h=h(R)$ so that 
a logarithmic fraction of $\{|\psi_R f|\geq\d^{\e_3}\l\}$ 
is contained in the set $\{|hg^R_h|\geq\d^{2\e_3}\l\}$. 
Finally, we pigeonhole to get a value of $h$ so that
a logarithmic fraction of $\{|f|\geq \l\}$ is contained in
$\bigcup_R \{|hg^R_h|\geq\d^{2\e_3}\l\})$.

Let $\l_*=\d^{2\e_3}\l h^{-1}$ and $f_R=g^R_h$,
with this value of $h$.  The lower bound in \bref{vp.lambda}
follows from the bound on $h$ just stated; for the upper bound, we use that
\[
\l_*=\d^{2\e_3}\l h^{-1}\leq \|g_h^R\|_{\infty}
\lesssim M_{\sigma,\rho}.
\]

Let $\P_R$ be the plate family for $f_R$.
From \bref{wa1} we have
\begin{equation}
\sum_{\pi\in\P_R}|\pi|\lesssim h^{-p} \|\psi_R f\|_{p,\rho}^p.
\label{ls14}\end{equation}
Letting $p=2$ and recalling the definition of $\l_*$, 
we deduce \bref{vp2}.

It remains to prove \bref{vp3}.  By \bref{ls14}, it suffices to show that
\begin{equation}\label{zzz.10}
\sum_R \|\psi_R f\|_{p,\rho}^p
\lesssim 
\d^{-C\e_3}M_{\rho,\d}^{rp+\alpha}\|f\|_2^2.
\end{equation}

We first claim that
\begin{equation}\label{ls10}
\sum_R \|\psi_R f \|_{p, \rho}^p
\lesssim\|f\|_{p,\rho}=\sum_{b}\|\Psi_b\ast f\|_p^p,
\end{equation}
where $\{\Psi_b\}$ is the partition of unity 
defining $\|\cdot\|_{p,\rho}$.
Indeed, observe that $\psi_R(\Xi_{a}*f)$ is Fourier supported in 
$\Pi_a+\calr$.  As shown in the proof of Lemma \ref{box-convolute},
the latter set is contained in a $C\rho$-neighbourhood of $\Pi_a$.
Hence the number of $a'$ such that
$ \Xi_{a'}*(\psi_R \cdot (\Xi_{a}*f))\neq 0$
is bounded by a constant independent of $\d$ and
$a$, and similarly with $a$ and $a'$ interchanged.
We thus have
\begin{eqnarray*}
\sum_R \|\psi_Rf\|_{p,\d}^p &=&
\sum_R\sum_a \|\Xi_a*(\psi_R f)\|_p^p
\\
&\lesssim&\sum_b\sum_R \sum_{b'}
\|\Xi_a\ast(\psi_R \cdot (\Xi_{a'}\ast f))\|_p^p
\\
&\lesssim&\sum_R\sum_{a'}\|\psi_R \cdot(\Xi_{a'}\ast f)\|_p^p
\\
&\lesssim&\sum_{a'}\|\Xi_{b'}\ast f\|_p^p,
\end{eqnarray*}
as required.

Note that $\supp\widehat{\Psi_b\ast f}\subset S^\rho_{C\d}$.  We also have
\[
\|\Psi_b\ast f\|_{\infty, \d}
=\max_a\|\Xi_a\ast \Psi_b\ast f\|_{\infty}
\lesssim\max_a\|\Psi_b\|_1\|\Xi_a \ast f\|_{\infty}
\lesssim\|f\|_{\infty, \delta}\lesssim 1.
\]
Applying the inductive assumption $P(p,\alpha)$ in the 
form \bref{mnmnmn} to $\Psi_b\ast f$, we see that
\[
\|\Psi_b\ast f\|_p^p\lesssim
\d^{-\e-C\e_3}M_{\rho,\d}^{rp+\alpha}\|\Psi_b\ast f\|_2^2.
\]
Combining this with \bref{ls10} and using the essential orthogonality
of $\Psi_b * f$, we obtain \bref{zzz.10} as claimed.
\hfill$\square$

\begin{lemma}\label{Lemma 5.2}
Assume that $P(p,\alpha,\e)$ holds, and that $f\in\ssd$  localizes
at $\l$.  Let $\P$ be the plate family for $f$.  
Then for any $\beta>(1-\e_0)\alpha$ we have
\begin{equation}\label{4.2}
|\{ |f|>\l\}|\lesssim \l^{-p} \d^{-\epsilon}
M_{\sigma,\d}^{rp+\beta}\sum_{\pi\in\P}|\pi|.
\end{equation}
\end{lemma}

\underline{Proof} Let $W=\{|f|\geq\l\}$.  The localization
assumption means that $f$ has subfunctions $f^Q$, where $Q$
ranges over $t\d^{-1}$-cubes, such that \bref{vo9} holds and
\[
|W|\lessapprox |\bigcup_Q W_Q|,
\]
where $W_Q=Q\cap\{|f_Q|\gtrapprox \l\}$.

Let $g_Q=\psi_Qf^Q$.  By Lemma \ref{Lemma 4.2}, $\supp g_Q\subset
S_{\d/t}$ and
$\|g_Q\|_{\infty, \d/t}\lesssim M_{\d/t,\d}$.  Applying 
the inductive hypothesis \bref{4.1} to $M_{\d/t,\d}^{-1}g_Q$,
with $\d$ replaced by $\d/t$ and $\l$ replaced by $\lo{-C}
M_{\d/t,\d}^{-1}\l$, we obtain 
\begin{eqnarray*}
|\{|g_Q|\gtrapprox \l\}|
&\lessapprox &(M_{\d/t,\d}^{-1}\l)^{-p}(\d/t)^{-\epsilon}
M_{\sigma,\d/t}^{rp+\alpha}\|M_{\d/t,\d}^{-1}g_Q\|_2^2 
\\
&\lessapprox&\l^{-p} M_{\d/t,\d}^{p-2}M_{\sigma,\d/t}^{rp+\alpha} 
(\d/t)^{-\epsilon}\|g_Q\|_2^2 .
\end{eqnarray*}
By \bref{Schur}, \bref{ls3}, and \bref{vo9}, we have
\[
\sum_Q \|g_Q\|_2^2\lesssim t\sum_Q \|f^Q\|_2^2
\lesssim t\sum_Q \sum_{\pi\in\P(f^Q)}|\pi|
\lessapprox t \sum_{\pi\in\P}|\pi|.
\]
Hence
\begin{eqnarray*}
|W|\lessapprox \sum_Q|W_Q|
&\lessapprox & \sum_Q|\{|g_Q|\gtrapprox \l\}|\\
&\lessapprox&
\l^{-p}(\d/t)^{-\epsilon} M_{\d/t,\d}^{p-2}M_{\sigma,\d/t}^{rp+\alpha}
t\sum_{\pi\in\P}|\pi|
\\
&\lessapprox & \l^{-p}(\d/t)^{-\epsilon}
M_{\sigma,\d}^{rp+\alpha}
\sum_{\pi\in\P}|\pi|
\cdot t M_{\d/t,\d}^{p-2-rp-\alpha},
\end{eqnarray*}
where at the last step we used that $M_{\d/t,\d}M_{\sigma,\d}=M_\d$.
Recall that $p-2-rp=\frac{2}{d-k}$ and $t=(\d/\sigma)^{\e_0}$.   
Using also \bref{eq-M}, we get
$$
t M_{\d/t,\d}^{p-2-rp-\alpha}=
t M_{\d/t,\d}^{\frac{2}{d-k}-\alpha}
\approx t \cdot (t^{-\frac{d-k}{2}})^{\frac{2}{d-k}-\alpha}$$
$$\approx  t^{(d-k)\alpha/2}
\approx  \big((\frac{\sigma}{\delta})^{-(d-k)/2}\big)^{-\e_0\alpha}
\approx M_{\sigma,\d}^{-\epsilon_0\alpha},$$
which yields \bref{4.2} as required.
\hfill$\square$

\begin{lemma}\label{Lemma 5.3}
Assume that we have $P(p,\alpha,\e)$, and that $f\in\ssd$
with plate family $\P$ satisfies
\begin{equation}\label{vo5}
|\P|\lesssim t^C\l^4(\d/\sigma)^{\frac{3d+1}{4}-k}.
\end{equation}
Then for any $\gamma>(1-\frac{\e_0}{2})\alpha$ we have
\begin{equation}\label{4.2a}
|\{ |f|>\l\}|\lesssim \l^{-p} \d^{-2\e-C\e_3}
M_{\sigma,\d}^{rp+\gamma}\sum_{\pi\in\P}|\pi|.
\end{equation}
\end{lemma}

\underline{Proof}
Apply Lemma \ref{Lemma 3.2} to $f$.  If $f$ localizes at $\l$, then we are done
by Lemma \ref{Lemma 5.2}.  Otherwise, pick a subfunction $f^*$ 
as in Lemma \ref{Lemma 3.2}, and apply Lemma \ref{Lemma 5.1} to it.
We obtain $f_{R}\in\Sigma^\sigma_\rho$ and a value of $\lambda_*$ as in
\bref{vp.lambda} so that
$|\{|f^*|\gtrapprox \l\}|\lessapprox |\bigcup_{R}
\{|f_{R}|\geq\l_*\})|$ and
\begin{equation}\label{zzz.20}
\sum_{\pi\in\P_R}|\pi|\leq\d^{-C\e_3}(\frac{\l_*}{\l})^2\|\psi_R f^*\|_2^2,
\end{equation}
\begin{equation}\label{zzz.21}
\sum_{R}\sum_{\pi\in\P_R}|\pi|\leq\d^{-\e-C\e_3}
(\frac{\l_*}{\l})^{p}
M_{\rho,\delta}^{rp+\alpha}\|f^*\|_2^2
\end{equation}
From Lemma \ref{Lemma 3.2} and \bref{zzz.20} we have
\[
|\P_R|\rho^{-\frac{d-k}{2}-1}
\approx \sum_{\pi\in\P_R}|\pi|
\lesssim\d^{-C\e_3}\frac{\l_*^2}{\l^4}(\frac{\sigma}{\d})^{-k/2}
M_{\sigma,\d}|R||\P|.
\]
Plugging in \bref{eq-Msd} and \bref{R-eq2}, 
we obtain after some algebra that
\[
|\P_R|
\lesssim\d^{-C\e_3}\frac{\l_*^2}{\l^4}(\frac{\sigma}{\d})^{\frac{3d+1}{4}-k}
|\P|.
\]
Combining this with \bref{vo5}, we see that
\[
|\P_R|\leq\d^{-C\e_3}t^{C}\l_*^2.
\]
By Lemma \ref{Lemma 3.1} and \bref{sd-size}, $f_R$ localize at $\l_*$.  Applying 
Lemma \ref{Lemma 5.2} we see that
\begin{equation}\label{wb1}
|\{|f_R|\geq\l_*\}| \lesssim \d^{-\e}\l_*^{-p} 
M_{\sigma,\rho}^{rp+\beta}\sum_{\pi\in\P_R}|\pi|
\end{equation}
for any $\beta$ with $\beta>(1-\e_0)\alpha$.
Hence
\begin{eqnarray*}
|\{|f|\geq\l\}|&\lessapprox&
\sum_R |\{|f_R|\geq\l_*\}|
\\
&\leq&\d^{-2\e-C\e_3} \l_*^{-p}
M_{\sigma,\rho}^{rp+\beta}\sum_R\sum_{\pi\in\P_R}|\pi|
\\
&\leq&\d^{-2\e-C\e_3}\l^{-p}
M_{\sigma,\rho}^{rp+\beta}
M_{\rho,\delta}^{rp+\alpha}\|f^*\|_2^2
\\
&\leq&\d^{-2\e-C\e_3}\l^{-p}
M_{\sigma,\d}^{rp+(\alpha+\beta)/2}
\sum_{\pi\in\P}|\pi|,
\end{eqnarray*}
where we also used \bref{zzz.21} and Lemma 4.1 with $p=2$. 
\hfill$\square$

\bigskip

\underline{Proof of Proposition \ref{scales}.}
Assume that $\supp\widehat{f}\subset S^\sigma_\delta$ and 
\begin{equation}\label{pink-e1}
\|f\|_{\infty, \d}\leq 1.
\end{equation}
We observed at the beginning of the proof of Proposition \ref{cor-iterate}
that \bref{4.1} follows from Chebyshev's inequality if
\bref{cheb} holds.  In particular, it holds for any $\alpha,\e>0$ if
$\l\lesssim M_{\sigma,\d}^{rp/(p-2)}$.  We may therefore assume that
\begin{equation}
\l\geq M_{\sigma,\d}^{rp/(p-2)}.
\label{4.3}\end{equation}
Recall also that the assumptions of Proposition \ref{scales} include 
\bref{sd-size}.

Choose $f_R\in\Sigma^\sigma_\rho$ as in in Lemma \ref{Lemma 5.1},
with $\rho=\sqrt{\sigma\d}$. 
Suppose that we can prove that $f_R$ obey \bref{4.2a}, with $\d$ and $\l$
replaced by $\rho$ and $\l_*$.  Then, using also \bref{vp3} and
\bref{eq-Msd}, we obtain
\begin{eqnarray*}
|\{|f|\geq\l\}|&\lessapprox &
\sum_R |\{|f_R|\geq\l_*\}|
\\
&\lesssim &\l_*^{-p}\d^{-2\e-C\e_3}M_{\sigma,\rho}^{rp+\gamma}
\sum_{R}\sum_{\pi\in\P_{R}}|\pi|
\\
&\lesssim &\l^{-p}\d^{-3\e-C\e_3}M_{\sigma,\rho}^{rp+\gamma}
M_{\rho,\delta}^{rp+\alpha}\|f\|_2^2
\\
&\lesssim  &\l^{-p}\d^{-3\e-C\e_3}M_{\sigma,\delta}^{rp+\theta}\|f\|_2^2
\end{eqnarray*}
for any $\theta>(\alpha+\gamma)/2$, as required.

It thus remain to find conditions under which $f_R$ obey the assumptions
of either Lemma \ref{Lemma 5.2} or \ref{Lemma 5.3}.
From \bref{vp2}, \bref{ls30}, \bref{pink-e1} we have
\begin{eqnarray*}
|\P_{R}|\rho^{-\frac{d-k}{2}-1}
&\approx& \sum_{\pi\in\P_R}|\pi|
\\
&\lesssim& \d^{-C\e_3}(\frac{\l_*}{\l})^2 \|\psi_R f\|_2^2
\\
&\lesssim& \d^{-C\e_3}(\frac{\l_*}{\l})^2 M_{\sigma,\d} |R|.
\end{eqnarray*}
Plugging in the expression \bref{R-eq2} for $|R|$, we obtain 
after minor simplifications
\begin{equation}\label{zzz.100}
|\P_{R}|\lesssim 
\d^{-C\e_3}(\frac{\l_*}{\l})^2 M_{\sigma,\d}(\frac{\sigma}{\d})^{\frac{d+1}{4}}. 
\end{equation}

Assume first that $d> 3k$ and 
$p>p_1(d,k):=2+\frac{8}{d-3k}$.  From \bref{zzz.100} 
and \bref{4.3} we have
\[
|\P_{R}|
\lesssim \d^{-C\e_3}\l_*^2 M_{\sigma,\d}^{1-\frac{2rp}{p-2}}
(\frac{\sigma}{\d})^{\frac{d+1}{4}}. 
\]
Recalling \bref{def-r}, we see that $1-\frac{2rp}{p-2}=
-1+\frac{4}{(d-k)(p-2)}$.
We now also plug in \bref{eq-Msd} for $M_{\sigma,\d}$.
After some algebra, this yields 
\[
|\P_R|\lesssim 
\d^{-C\e_3}\l_*^2 
(\frac{\sigma}{\d})^{-\frac{d-3k}{4}+\frac{2}{p-2}}.
\]
Recall that $t=(\d/\sigma)^{\e_0}$.  
If we assume that $\e_0$ is sufficiently small depending on $p,d,k$, and
if $\e_3$ is small enough compared to the $\e_2$ in \bref{sd-size}, then
our assumption that $p>p_1(d,k)$ implies that
$|\P_R|\lesssim t^C \l_*^2$.  By Lemma \ref{Lemma 3.1}
$f_R$ localize, hence Lemma \ref{Lemma 5.2} applies.

Assume now that $p>p_2(d,k):=2+\frac{32}{3d-4k-3}$; in this case
we will see that \bref{vo5} holds, with $\lambda$ and $\d$ replaced by
$\l_*$ and $\rho$.  It suffices to check that the
right side of \bref{zzz.100} is 
$\lesssim t^C\l_*^4(\rho/\sigma)^{\frac{3d+1}{4}-k}$.
After some simplifications, this is equivalent to
\[
\l_*^2\l^2\gtrsim \d^{-C\e_3}t^{-C}M_{\sigma,\d}(\frac{\sigma}{\d})^{\frac{5d+3-4k}{8}}.
\]
By \bref{4.3} and \bref{vp.lambda}, we have
\[
\l^2\l_*^2\gtrsim \d^{C\e_3}\l^4 M_{\rho,\d}^{-2}
\gtrsim \d^{C\e_3}M_{\sigma,\d}^{4rp/(p-2)} M_{\rho,\d}^{-2}
\approx \d^{C\e_3}M_{\sigma,\d}^{\frac{4rp}{p-2}-1},
\]
where at the last step we used \bref{eq-Msd}.  
It thus suffices to prove that
\[
M_{\sigma,\d}^{\frac{4rp}{p-2}-1}
\gtrsim \d^{-C\e_3}t^{-C}M_{\sigma,\d}(\frac{\sigma}{\d})^{\frac{5d+3-4k}{8}}.
\]
Using again \bref{eq-Msd} and \bref{def-r}, we see after some algebra
that this is equivalent to
\[
(\frac{\sigma}{\d})^{\frac{3d-4k-3}{8}-\frac{4}{p-2}}
\gtrsim \d^{-C\e_3}t^{-C}.
\]
But if $p>p_2(d,k)$, then the exponent on the left is positive, hence the
last estimate follows if $\e_0$ was chosen small enough and if $\e_3$ is
small enough compared to $\e_2$ in \bref{sd-size}.
\qed


\section{Properties of k-cones}\label{k-cones}


We first recall the construction of $k$-cones described in 
the introduction.   In this section, we will use superscripts 
to denote Cartesian coordinates in $\rr^{d+1}$, e.g.
$x=(x^1,\dots,x^{d+1})$.  

Let $L_0$ be a $(d-k+1)$-dimensional linear subspace of $\rr^{d+1}$,
and let $L_i=L_0+v_i$ for $i=1,\dots,k$, where $v_1,\dots,v_{k}$ are
linearly independent vectors such that $L_0,v_1,\dots,v_k$ span
$\rr^{d+1}$.  Applying an affine transformation if necessary,
we may assume that 
$$L_0=\{(0,\dots,0,x^{k+1},\dots,x^{d+1}):\ x^{k+1},\dots,x^{d+1}\in\rr\},$$
and that for $i=1,\dots,k$, 
$$v_i=(v_i^1,\dots,v_i^{d+1}),\ v_i^i=1,\ v_i^j=0\hbox{ if }i\neq j.$$
If $x=(x^1,\dots,x^{d+1})\in\rr^{d+1}$, 
we will also use the notation $x=(x^\perp,x^\parallel)$,
where $x^\perp=(x^1,\dots,x^k)$ and $x^\parallel=(x^{k+1},\dots,x^{d+1})$
denote the components orthogonal to $L_0$ and parallel to it, respectively.

We let $E_i\subset L_i$, $i=0,\dots,k$, be surfaces of dimension
$d-k$ such that $E_i$ is the boundary of a strictly convex solid
in $L_i$, is smooth, and has nonvanishing Gaussian curvature.  Thus for each unit vector
$n\in S^{d-k}$, each $E_i$ contains exactly one point $x_i$ such
that $n$ is the outward unit normal vector to $E_i$ in $L_i$ at $x_i$.
We will then write $n=n_i(x_i)$.  Since $E_i$ is smooth, the mapping
$x_i\to n_i(x_i)$ is a smooth diffeomorphism from $E_i$ to $S_i$.

We say that a $(k+1)$-tuple of points $(x_0,x_1,\dots,x_k)$, 
$x_i\in\rr^{d+1}$, is {\em good} if $x_i\in E_i$, $i=0,\dots,k$, and if 
the outward unit normal vectors to $E_i$ in $L_i$ at $x_i$ are the same
(i.e. $n_0(x_0)=\dots=n_k(x_k)$).
We then let
$$
S=\bigcup_{(x_0,\dots,x_k)\ good}\eta(x_0,\dots,x_k),
$$
where $\eta(x_0,\dots,x_k)$ denotes the $k$-dimensional convex hull
of $x_0,\dots, x_k$ in $\rr^{d+1}$.  


We first verify that $\eta(x_0,\dots,x_k)$ is indeed $k$-dimensional.
Indeed, if $\eta(x_0,\dots,x_k)$ had dimension less than $k$, then 
the dimension of the affine space spanned by 
$\eta(x_0,\dots,x_k)$ and $L_0$ would be less than $d+1$.  But on 
the other hand, this affine space contains both $L_0$ and $v_1,
\dots,v_k$ (since $v_i-x_i\in L_i$), hence must be equal to all
of $\rr^{d+1}$, which proves our claim.  

Let $T_i(x_i)$ be the $d-k$-dimensional affine space tangent to
$E_i$ in $L_i$ at $x_i$.  Note that if $(x_0,\dots,x_k)$ are good,
then $T_i(x_i)$, $i=0,\dots,k$, are parallel.
The above argument shows that $T_i(x_i), n_i(x_i), L_i$ span $\rr^{d+1}$.
In fact, we can say slightly more.  
Let $D$ be a closed disc in $S^{d-k}$, and let 
$$
E_{i}|_D=\{x_i\in E_i:\ n_i(x_i)\in D\},
$$
$$
S|_D=\bigcup_{(x_0,\dots,x_k)\ good,\ n_i(x_i)\in D}\eta(x_0,\dots,x_k).
$$
We may then choose orthonormal bases 
\begin{equation}\label{f-e1}
F^\parallel(x_0)=\{u_1(x_0),\dots,u_{d-k}(x_0)\},
\end{equation}
\begin{equation}\label{f-e2}
F^\perp(x_0)=\{u_{d-k+1}(x_0),\dots,u_{d}(x_0)\} ,
\end{equation}
for $T_0(x_0)$ and $\eta(x_0,\dots,x_k)$, respectively,  so that
each $u_j(x_0)$ depends smoothly on $x_0\in E_0|_D$.
(For $F^\parallel$, this is clear from smoothness of $S$; for
$F^\perp$, it can be done by applying Gram-Schmidt orthonormalization
to $x_1-x_0,\dots,x_k-x_0$.)
Let $V(x_0)$ denote the volume of a parallelepiped
spanned by $F^\parallel(x_0)$, $F^\perp(x_0)$, $n_0(x_0)$.
From the above considerations we have $V(x_0)\neq 0$.  But also
$V(x_0)$ is a continuous function of $x_0$, hence
\begin{equation}\label{nondeg}
V(x_0)\geq c_0>0,\ x_0\in E_0|_D.
\end{equation}

As noted above, $\eta(x_0,\dots,x_k)$ depend smoothly on 
$(x_0,\dots,x_k)$.  Thus to show that $S$ is a smooth surface of
codimension 1 in $\rr^{d+1}$, it suffices to prove that it is
nonsingular, i.e. the different $\eta(x_0,\dots,x_k)$ do not intersect.
This follows from the lemma below.

\begin{lemma}\label{nointersect}
Let $L$ be a $(d-k+1)$-dimensional affine subspace parallel to $L_0$.
If $S$ intersects $L$, then the cross-section $E:=S\cap L$ is 
a closed smooth $d-k$-dimensional surface with nonvanishing
Gaussian curvature, bounding a strictly convex body in $L$.
Moreover, each $\eta(x_0,\dots,x_k)$ intersects $L$ at
a unique point $x\in E\cap\eta(x_0,\dots,x_k)$.  The mapping $x_0\to x$
is smooth, and the unit
outward normal vector to $E$ in $L$ at $x$ is equal to $n_0(x_0)$.
\end{lemma}

\underline{Proof} 
Suppose that 
$z=(z^1,\dots,z^{d+1})\in L\cap \eta(z_0,\dots,z_k)$, and that
$z\notin\{z_0,\dots,z_k\}$.
Since $z$ belongs to the convex hull of $z_0,\dots,z_k$, taking
the first $k$ coordinates we may write
$$z^\perp=\alpha_0(0,\dots,0)+\alpha_1(1,0,\dots,0)
+\dots+\alpha_k(0,\dots,0,1),$$
with all $\alpha_i\geq 0$, $\sum_{j=0}^k\alpha_j=1$,
and at least one $\alpha_j>0$.  Since along
$\eta(z_0,\dots,z_k)$ each $z^j$, $j>k$, is a linear function of
$z^1,\dots,z^k$, we must in fact have
$$z=\alpha_0z_0+\dots+\alpha_kz_k.$$
Similarly, if $x$ is any other point in $L\cal S$ with
$x\in L\cap \eta(x_0,\dots,x_k)$, then $x^\perp=z^\perp$, so that
$$x=\alpha_0x_0+\dots+\alpha_kx_k,$$
with the same $\alpha_0,\dots,\alpha_k$.
This shows that $x$ depends smoothly on $x_0$.  

Next, we show that the mapping $x_0\to x$ is one-to-one.  Suppose
to the contrary that $x\in L\cap\eta(x_0,\dots,x_k)\cap\eta(y_0,
\dots,y_k)$.  Then
\begin{equation}\label{k-extreme}
x^\parallel=
\alpha_0x_0^{\parallel}+\dots +\alpha_kx_k^{\parallel}
=\alpha_0y_0^{\parallel}+\dots +\alpha_ky_k^{\parallel}.
\end{equation}
Consider, however, the scalar product 
\begin{equation}\label{k-extreme2}
n_0(x_0)\cdot(\alpha_0w_0^{\parallel}+\dots +\alpha_kw_k^{\parallel}),
\end{equation}
where $(w_0,\dots,w_k)$ ranges over all good
$k$-tuples.  By the strict convexity of each $E_i$, \bref{k-extreme2}
is maximized when $(w_0,\dots,w_k)=(x_0,\dots,x_k)$, and only
there.  But this contradicts the second equality in \bref{k-extreme}.

Interpreting \bref{k-extreme2} as the distance from the point
$w\in L\cap\eta(w_0,\dots,w_k)$ to the plane $\{y\in L: n_0(x_0)\cdot 
(y^\parallel-x^\parallel)=0$, and observing that it is a smooth
function of $w_0$, we see that $E$ is indeed a smooth surface in
the neighbourhood of $x_0$.
The argument in the last paragraph now shows that $n_0(x_0)$ is
indeed the unit outward normal vector to $E$ in $L$ at $x$.
\qed

\medskip

We now define a covering of $S_\d$, $\d>0$, by $\d$-sectors as follows.
We may assume that $S$ is contained in a large cube of sidelength $C_0$.
Fix a $\d^{1/2}$-separated subset ${\cal N}$ of the unit sphere $S^{d-k}$.
For $i=0,1,\dots,k$, let $\M_i$ be the set of points in $E_i$  
where the outward unit normal vectors in $L_i$ belong to ${\cal N}$.
We then define $\M$ to be the set of the centers of mass of 
$\eta(x_0,\dots,x_k)$, $\ x_i\in\M_i$.
We also define the corresponding $\d$-sectors to be rectangular 
boxes $\Pi_a$ centered at $a\in\M$ such that if $a$ is the center of mass
of $\eta(x_0,\dots,x_k)$, then $\Pi_a$ has sidelengths $C\times\dots\times C$
in the directions parallel to $\eta(x_0,\dots,x_k)$, 
$C\d^{1/2}\times\dots\times C\d^{1/2}$ in directions
tangent to $S$ at $a$ but perpendicular to $\eta(x_0,\dots,x_k)$, 
and $C\d$ in the direction parallel to $n(a)$, the normal vector
to $S$ at $a$.  Here $C$ is a large constant with $C>2C_0$, to be fixed later.

First, let $\Pi_a^0=\Pi_a\cap L_0$.  Then $\Pi_a^0$ is 
a rectangular box, centered at $x_0$, of dimensions at least $C\d^{1/2}\times\dots
\times C\d^{1/2}\times C\d$, where the long axes are tangent
to $E_0$ at $x_0$.  We claim that the dimensions of this box cannot exceed
$C'\d^{1/2}\times\dots \times C'\d^{1/2}\times C'\d$, for some other
constant $C'$.  To this end it suffices to prove that the maximal
number of disjoint translates of the box spanned by $(\d n_0(x_0),
\d u_1(x_0),\dots,\d u_{d-k}(x_0))$ (recall that $u_i$ were defined in
\bref{f-e1}) that can be placed inside $\Pi_a^0$ is bounded by
$\lesssim\d^{-(d-k)/2}$.  But consider the corresponding translates
of the box spanned by $(\d n_0(x_0),
\d u_1(x_0),\dots,\d u_{d-k}(x_0), u_{d-k+1}(x_0),\dots, u_{d}(x_0))$.
They are also disjoint, each one has volume $\gtrsim \d^{d-k+1}$ (by
\bref{nondeg}), and they are all contained in $\Pi_a$ which has
volume $\approx \d^{\frac{d-k}{2}+1}$.
Thus the claim follows by comparing volumes.

This shows that $\{\Pi_a^0\}_{a\in\M_0}$ is, if $C$ is large enough,
a standard finitely overlapping covering of the $\d$-neighbourhood of
$E_0$ in $L_0$ by rectangular boxes of dimensions roughly $C\d^{1/2}\times\dots
\times C\d^{1/2}\times C\d$, with centers in a $\d^{1/2}$-separated
set $\M_0$.  The same argument can now be repeated for $E_1,\dots,
E_k$.  In particular, this implies the finite overlap property.

It remains to prove the angular separation property.  We may restrict 
our attention to a small segment $S|_D$ as defined earlier such
that $n_0(x_0)\neq -n_0(y_0)$ for any $x_0,y_0\in E_0|_D$.  
We then want to prove that if $x,y\in\M$, $x\neq y$, then $|n(x)
-n(y)|\gtrsim\d^{1/2}$.  It suffices to prove that 
$$|n(x)^\parallel-n(y)^\parallel|\gtrsim \d^{1/2}.$$
But if $x\in\eta(x_0,\dots,x_k)$, $y\in\eta(y_0,\dots,y_k)$, 
then $n(x)^\parallel=c(x_0)n_0(x_0)$ (since both vectors are
orthogonal to $T_0(x_0)$), and similarly $n(y)^\parallel=c(y_0)n_0(y_0)$.
By the nonvanishing curvature assumption for $E_0$, the angle between 
$n_0(x_0)$ and $n_0(y_0)$ is at least $\d^{1/2}$; hence it suffices
to prove that $c(x_0),c(y_0)\gtrsim 1$.  But on the other hand,
we have $c(x_0)=V(x_0)$, where $V(x_0)$ is the volume of the 
parallelepiped defined before \bref{nondeg}, and similarly for
$y_0$.  Thus the claim follows from \bref{nondeg}.


\noindent{\sc Department of Mathematics, University of British Columbia, Vancouver, 
B.C. V6T 1Z2, Canada}

\noindent{\it ilaba@math.ubc.ca}

\bigskip

\noindent{\sc Department of Mathematics, University of Rochester,
Rochester, NY 14627, U.S.A.}

\noindent{\it malabika@math.rochester.edu}

\end{document}